\documentclass[11pt]{article}
\usepackage{xypic}
\usepackage{mathrsfs}
\usepackage{amssymb}
\usepackage{amsfonts}
\usepackage{amsmath}

\def\ben{\begin{enumerate}}
\def\een{\end{enumerate}}

\def\bit{\begin{itemize}}
\def\eit{\end{itemize}}

\def\0{\leqno}

\input xy
\xyoption{all}
\input{tcilatex}
\begin{document}

\begin{center}
{\Large {THE METRIZABILITY OF THE GENERALIZED TANGENT }\ \ \bigskip }

{\Large {BUNDLE OF DUAL OF A VECTOR BUNDLE }\ \ \bigskip }

\textbf{by }

\textbf{CONSTANTIN M. ARCU\c{S} }
\end{center}

\bigskip

\bigskip





%

\ \

\begin{abstract}
Two new classes of metrizable vector bundles have been presented in the
papers $[1]$ and $[4]$. The Lie algebroid generalized tangent bundle of a
dual vector bundle is presented. This Lie algebroid is a new example of
metrizable vector bundle. A new class of Hamilton spaces, called by use,
generalized Hamilton $\left( \rho ,\eta \right) $-space, Hamilton $\left(
\rho ,\eta \right) $-space and Cartan $\left( \rho ,\eta \right) $-space are
presented. The results obtained in the particular case of Lie algebroids
emphasize the importance and the utility of our new method by work. In
particular, if all morphisms are identities morphisms, then the classical
results are obtained. \ \ \bigskip\newline
\end{abstract}

\tableofcontents

\bigskip \bigskip \bigskip

\bigskip \bigskip \bigskip

\bigskip \bigskip \bigskip

\section{Introduction}

The study of the geometry of the usual Lie algebroid
\begin{equation*}
\left( \left( TT^{\ast }M,\tau _{T^{\ast }M},T^{\ast }M\right) ,\left[ ,%
\right] _{TT^{\ast }M},\left( Id_{TT^{\ast }M},Id_{T^{\ast }M}\right) \right)
\end{equation*}%
with a metrical structure
\begin{equation*}
\begin{array}{c}
g=g^{ij}dp_{i}\otimes dp_{j}\in \mathcal{T}~_{0}^{2}\left( VTT^{\ast }M,{%
\tau _{T^{\ast }M}},T^{\ast }M\right) ,%
\end{array}%
\end{equation*}%
was extensively examined by geometers and physicists in the framework of
generalized Hamilton space. (see $\left[ 15\right] $).

We know that a regular Hamiltonian on $T^{\ast }M$ is a smooth function $%
\begin{array}[b]{ccc}
T^{\ast }M & ^{\underrightarrow{~H\ }} & \mathbb{R}%
\end{array}%
$ such that the Hessian matrix with entries%
\begin{equation*}
\begin{array}[b]{c}
g^{ij}=\frac{1}{2}\frac{\partial ^{2}H}{\partial p_{i}\partial p_{j}}%
\end{array}%
\end{equation*}%
is everywhere nondegenerate. If the metrical structure of a generalized
Hamilton space is determined by a regular Hamiltonian, then we obtain the
Hamilton space. The concept of Hamilton, introduced by R. Miron in $\left[
14,13\right] ,$ vas intesinvely studied in $\left[ 6,7,8,10,16,17,...\right]
$ and it has been succesful as a geometric theory of the Hamiltonian
fundamental function, the fundamental entity in Mechanics and Physics$.$ In
the general framework of generalized Lie algebroids, the geometry of the
Hamilton fundamental function has been developed in the paper $\left[ 3%
\right] .$

If $H$ is square of a function on $T^{\ast }M,$ positively, $1$-homogeneous
with respect to the momentum $p_{i},$ then an important class of Hamilton
spaces, called by use Cartan spaces, were introduced by R. Miron $\left[
11,12\right] $. The geometry of Cartan space is a subgeometry of the
geometry of the Lie algebroid
\begin{equation*}
\left( \left( TT^{\ast }M,\tau _{T^{\ast }M},T^{\ast }M\right) ,\left[ ,%
\right] _{TT^{\ast }M},\left( Id_{TT^{\ast }M},Id_{T^{\ast }M}\right)
\right) .
\end{equation*}

Important contributions to the geometry of Cartan spaces were obtained E.
Cartan $\left[ 5\right] $ and A. Kawaguchi $\left[ 9\right] ,...$

The study of the metrizability in the general framework of generalized Lie
algebroids was extensively studied in the papers $\left[ 1,4\right] $. Using
a generalized Lie algebroid, we obtain the Lie algebroid generalized tangent
bundle
\begin{equation*}
\left( \left( \left( \rho ,\eta \right) T\overset{\ast }{E},\left( \rho
,\eta \right) \tau _{\overset{\ast }{E}},\overset{\ast }{E}\right) ,\left[ ,%
\right] _{\left( \rho ,\eta \right) T\overset{\ast }{E}},\left( \overset{%
\ast }{\tilde{\rho}},Id_{\overset{\ast }{E}}\right) \right)
\end{equation*}%
of dual vector bundle $\left( \overset{\ast }{E},\overset{\ast }{\pi }%
,M\right) .$ Using the basic notions and results presented in Sections $2,$ $%
3$ and $4$ we study the metrizability of this Lie algebroid in Section $5$.
In the particular case of Lie algebroids, we obtain important results.
Moreover, we obtain new results for the metrizability of the usual Lie
algebroid
\begin{equation*}
\left( \left( T\overset{\ast }{E},\tau _{\overset{\ast }{E}},\overset{\ast }{%
E}\right) ,\left[ ,\right] _{T\overset{\ast }{E}},\left( Id_{T\overset{\ast }%
{E}},Id_{\overset{\ast }{E}}\right) \right) .
\end{equation*}

Finally, in Section $6,$ we introduced a new class of Hamilton spaces,
called by use\emph{\ generalized Hamilton }$\left( \rho ,\eta \right) $\emph{%
-spaces, Hamilton }$\left( \rho ,\eta \right) $\emph{-spaces and Cartan }$%
\left( \rho ,\eta \right) $\emph{-spaces.}

In the particular case of Lie algebroids, new and important results are
obtained. In particular, if $\left( \rho ,\eta ,h\right) =\left(
Id_{TM},Id_{M},Id_{M}\right) ,$ then the classical results are obtained.

\section{Preliminaries}

Let$\mathbf{~Vect},$ $\mathbf{Liealg},~\mathbf{Mod}$\textbf{,} $\mathbf{Man}$
and $\mathbf{B}^{\mathbf{v}}$ be the category of real vector spaces, Lie
algebras, modules, manifolds and vector bundles respectively.

We know that if $\left( E,\pi ,M\right) \in \left\vert \mathbf{B}^{\mathbf{v}%
}\right\vert ,$ $\Gamma \left( E,\pi ,M\right) =\left\{ u\in \mathbf{Man}%
\left( M,E\right) :u\circ \pi =Id_{M}\right\} $ and $\mathcal{F}\left(
M\right) =\mathbf{Man}\left( M,\mathbb{R}\right) ,$ then $\left( \Gamma
\left( E,\pi ,M\right) ,+,\cdot \right) $ is a $\mathcal{F}\left( M\right) $%
-module. If \ $\left( \varphi ,\varphi _{0}\right) \in \mathbf{B}^{\mathbf{v}%
}\left( \left( E,\pi ,M\right) ,\left( E^{\prime },\pi ^{\prime },M^{\prime
}\right) \right) $ such that $\varphi _{0}\in Iso_{\mathbf{Man}}\left(
M,M^{\prime }\right) ,$ then, using the operation
\begin{equation*}
\begin{array}{ccc}
\mathcal{F}\left( M\right) \times \Gamma \left( E^{\prime },\pi ^{\prime
},M^{\prime }\right) & ^{\underrightarrow{~\ \ \cdot ~\ \ }} & \Gamma \left(
E^{\prime },\pi ^{\prime },M^{\prime }\right) \\
\left( f,u^{\prime }\right) & \longmapsto & f\circ \varphi _{0}^{-1}\cdot
u^{\prime }%
\end{array}%
\end{equation*}%
it results that $\left( \Gamma \left( E^{\prime },\pi ^{\prime },M^{\prime
}\right) ,+,\cdot \right) $ is a $\mathcal{F}\left( M\right) $-module and we
obtain the $\mathbf{Mod}$-morphism%
\begin{equation*}
\begin{array}{ccc}
\Gamma \left( E,\pi ,M\right) & ^{\underrightarrow{~\ \ \Gamma \left(
\varphi ,\varphi _{0}\right) ~\ \ }} & \Gamma \left( E^{\prime },\pi
^{\prime },M^{\prime }\right) \\
u & \longmapsto & \Gamma \left( \varphi ,\varphi _{0}\right) u%
\end{array}%
\end{equation*}%
defined by
\begin{equation*}
\begin{array}{c}
\Gamma \left( \varphi ,\varphi _{0}\right) u\left( y\right) =\varphi \left(
u_{\varphi _{0}^{-1}\left( y\right) }\right) ,%
\end{array}%
\end{equation*}%
for any $y\in M^{\prime }.$

Let $M,N\in \left\vert \mathbf{Man}\right\vert ,$ $h\in Iso_{\mathbf{Man}%
}\left( M,N\right) $ and $\eta \in Iso_{\mathbf{Man}}\left( N,M\right) $.

We know (see $\left[ 1,3,4\right] $) that if $\left( F,\nu ,N\right) \in
\left\vert \mathbf{B}^{\mathbf{v}}\right\vert $ so that there exists
\begin{equation*}
\begin{array}{c}
\left( \rho ,\eta \right) \in \mathbf{B}^{\mathbf{v}}\left( \left( F,\nu
,N\right) ,\left( TM,\tau _{M},M\right) \right)%
\end{array}%
\end{equation*}%
and also an operation
\begin{equation*}
\begin{array}{ccc}
\Gamma \left( F,\nu ,N\right) \times \Gamma \left( F,\nu ,N\right) & ^{%
\underrightarrow{\left[ ,\right] _{F,h}}} & \Gamma \left( F,\nu ,N\right) \\
\left( u,v\right) & \longmapsto & \left[ u,v\right] _{F,h}%
\end{array}%
\end{equation*}%
with the following properties:\bigskip

\noindent $\qquad GLA_{1}$. the equality holds good
\begin{equation*}
\begin{array}{c}
\left[ u,f\cdot v\right] _{F,h}=f\left[ u,v\right] _{F,h}+\Gamma \left(
Th\circ \rho ,h\circ \eta \right) \left( u\right) f\cdot v,%
\end{array}%
\end{equation*}%
\qquad \quad\ \ for all $u,v\in \Gamma \left( F,\nu ,N\right) $ and $f\in
\mathcal{F}\left( N\right) .$

\medskip $GLA_{2}$. the $4$-tuple $\left( \Gamma \left( F,\nu ,N\right)
,+,\cdot ,\left[ ,\right] _{F,h}\right) $ is a Lie $\mathcal{F}\left(
N\right) $-algebra,

$GLA_{3}$. the $\mathbf{Mod}$-morphism $\Gamma \left( Th\circ \rho ,h\circ
\eta \right) $ is a $\mathbf{LieAlg}$-morphism of
\begin{equation*}
\left( \Gamma \left( F,\nu ,N\right) ,+,\cdot ,\left[ ,\right] _{F,h}\right)
\end{equation*}%
source and
\begin{equation*}
\left( \Gamma \left( TN,\tau _{N},N\right) ,+,\cdot ,\left[ ,\right]
_{TN}\right)
\end{equation*}%
target, \medskip \noindent then the triple $\left( \left( F,\nu ,N\right) ,%
\left[ ,\right] _{F,h},\left( \rho ,\eta \right) \right) $ is called \emph{%
generalized Lie algebroid.\ }

In particular, if $h=Id_{M}=\eta ,$ then we obtain the definition of the Lie
algebroid.

Let $\left( \left( F,\nu ,N\right) ,\left[ ,\right] _{F,h},\left( \rho ,\eta
\right) \right) $ be an generalized Lie algebroid.

\begin{itemize}
\item Locally, for any $\alpha ,\beta \in \overline{1,p},$ we set $\left[
t_{\alpha },t_{\beta }\right] _{F,h}=L_{\alpha \beta }^{\gamma }t_{\gamma }.$
We easily obtain that $L_{\alpha \beta }^{\gamma }=-L_{\beta \alpha
}^{\gamma },~$for any $\alpha ,\beta ,\gamma \in \overline{1,p}.$
\end{itemize}

The real local functions $L_{\alpha \beta }^{\gamma },~\alpha ,\beta ,\gamma
\in \overline{1,p}$ will be called the \emph{structure functions of the
generalized Lie algebroid }$\left( \left( F,\nu ,N\right) ,\left[ ,\right]
_{F,h},\left( \rho ,\eta \right) \right) .$

\begin{itemize}
\item We assume the following diagrams:%
\begin{equation*}
\begin{array}[b]{ccccc}
F & ^{\underrightarrow{~\ \ \ \rho ~\ \ }} & TM & ^{\underrightarrow{~\ \ \
Th~\ \ }} & TN \\
~\downarrow \nu &  & ~\ \ \ \downarrow \tau _{M} &  & ~\ \ \ \downarrow \tau
_{N} \\
N & ^{\underrightarrow{~\ \ \ \eta ~\ \ }} & M & ^{\underrightarrow{~\ \ \
h~\ \ }} & N \\
&  &  &  &  \\
\left( \chi ^{\tilde{\imath}},z^{\alpha }\right) &  & \left(
x^{i},y^{i}\right) &  & \left( \chi ^{\tilde{\imath}},z^{\tilde{\imath}%
}\right)%
\end{array}%
\end{equation*}

where $i,\tilde{\imath}\in \overline{1,m}$ and $\alpha \in \overline{1,p}.$

If%
\begin{equation*}
\left( \chi ^{\tilde{\imath}},z^{\alpha }\right) \longrightarrow \left( \chi
^{\tilde{\imath}\prime }\left( \chi ^{\tilde{\imath}}\right) ,z^{\alpha
\prime }\left( \chi ^{\tilde{\imath}},z^{\alpha }\right) \right) ,
\end{equation*}%
\begin{equation*}
\left( x^{i},y^{i}\right) \longrightarrow \left( x^{i%
{\acute{}}%
}\left( x^{i}\right) ,y^{i%
{\acute{}}%
}\left( x^{i},y^{i}\right) \right)
\end{equation*}%
and
\begin{equation*}
\left( \chi ^{\tilde{\imath}},z^{\tilde{\imath}}\right) \longrightarrow
\left( \chi ^{\tilde{\imath}\prime }\left( \chi ^{\tilde{\imath}}\right) ,z^{%
\tilde{\imath}\prime }\left( \chi ^{\tilde{\imath}},z^{\tilde{\imath}%
}\right) \right) ,
\end{equation*}%
then
\begin{equation*}
\begin{array}[b]{c}
z^{\alpha
{\acute{}}%
}=\Lambda _{\alpha }^{\alpha
{\acute{}}%
}z^{\alpha }%
\end{array}%
,
\end{equation*}%
\begin{equation*}
\begin{array}[b]{c}
y^{i%
{\acute{}}%
}=\frac{\partial x^{i%
{\acute{}}%
}}{\partial x^{i}}y^{i}%
\end{array}%
\end{equation*}%
and
\begin{equation*}
\begin{array}{c}
z^{\tilde{\imath}\prime }=\frac{\partial \chi ^{\tilde{\imath}\prime }}{%
\partial \chi ^{\tilde{\imath}}}z^{\tilde{\imath}}.%
\end{array}%
\end{equation*}

\item We assume that $\left( \theta ,\mu \right) \overset{put}{=}\left(
Th\circ \rho ,h\circ \eta \right) $. If $z^{\alpha }t_{\alpha }\in \Gamma
\left( F,\nu ,N\right) $ is arbitrary, then
\begin{equation*}
\begin{array}[t]{l}
\begin{array}{c}
\Gamma \left( Th\circ \rho ,h\circ \eta \right) \left( z^{\alpha }t_{\alpha
}\right) f\left( h\circ \eta \left( \varkappa \right) \right) =\vspace*{1mm}
\\
=\left( \theta _{\alpha }^{\tilde{\imath}}z^{\alpha }\frac{\partial f}{%
\partial \varkappa ^{\tilde{\imath}}}\right) \left( h\circ \eta \left(
\varkappa \right) \right) =\left( \left( \rho _{\alpha }^{i}\circ h\right)
\left( z^{\alpha }\circ h\right) \frac{\partial f\circ h}{\partial x^{i}}%
\right) \left( \eta \left( \varkappa \right) \right) ,%
\end{array}%
\end{array}%
\leqno(2.1)
\end{equation*}%
for any $f\in \mathcal{F}\left( N\right) $ and $\varkappa \in N.$
\end{itemize}

The coefficients $\rho _{\alpha }^{i}$ respectively $\theta _{\alpha }^{%
\tilde{\imath}}$ change to $\rho _{\alpha
{\acute{}}%
}^{i%
{\acute{}}%
}$ respectively $\theta _{\alpha
{\acute{}}%
}^{\tilde{\imath}%
{\acute{}}%
}$ according to the rule:
\begin{equation*}
\begin{array}{c}
\rho _{\alpha
{\acute{}}%
}^{i%
{\acute{}}%
}=\Lambda _{\alpha
{\acute{}}%
}^{\alpha }\rho _{\alpha }^{i}\displaystyle\frac{\partial x^{i%
{\acute{}}%
}}{\partial x^{i}},%
\end{array}%
\leqno(2.2)
\end{equation*}%
respectively%
\begin{equation*}
\begin{array}{c}
\theta _{\alpha
{\acute{}}%
}^{\tilde{\imath}%
{\acute{}}%
}=\Lambda _{\alpha
{\acute{}}%
}^{\alpha }\theta _{\alpha }^{\tilde{\imath}}\displaystyle\frac{\partial
\varkappa ^{\tilde{\imath}%
{\acute{}}%
}}{\partial \varkappa ^{\tilde{\imath}}},%
\end{array}%
\leqno(2.3)
\end{equation*}%
where
\begin{equation*}
\left\Vert \Lambda _{\alpha
{\acute{}}%
}^{\alpha }\right\Vert =\left\Vert \Lambda _{\alpha }^{\alpha
{\acute{}}%
}\right\Vert ^{-1}.
\end{equation*}

\emph{Remark 2.1 } The following equalities hold good:%
\begin{equation*}
\begin{array}{c}
\displaystyle\rho _{\alpha }^{i}\circ h\frac{\partial f\circ h}{\partial
x^{i}}=\left( \theta _{\alpha }^{\tilde{\imath}}\frac{\partial f}{\partial
\varkappa ^{\tilde{\imath}}}\right) \circ h,\forall f\in \mathcal{F}\left(
N\right) .%
\end{array}%
\leqno(2.4)
\end{equation*}%
\emph{and }%
\begin{equation*}
\begin{array}{c}
\displaystyle\left( L_{\alpha \beta }^{\gamma }\circ h\right) \left( \rho
_{\gamma }^{k}\circ h\right) =\left( \rho _{\alpha }^{i}\circ h\right) \frac{%
\partial \left( \rho _{\beta }^{k}\circ h\right) }{\partial x^{i}}-\left(
\rho _{\beta }^{j}\circ h\right) \frac{\partial \left( \rho _{\alpha
}^{k}\circ h\right) }{\partial x^{j}}.%
\end{array}%
\leqno(2.5)
\end{equation*}

If $\left( E,\pi ,M\right) \in \left\vert \mathbf{B}^{\mathbf{v}}\right\vert
,$ then we have the $\mathbf{B}^{\mathbf{v}}$-morphism%
\begin{equation*}
\begin{array}{ccc}
~\ \ \ \ \ \ \ \ \ \ \ \ \ \overset{\ast }{\pi }^{\ast }\left( h^{\ast
}F\right) & \hookrightarrow & F \\
\overset{\ast }{\pi }^{\ast }\left( h^{\ast }\nu \right) \downarrow &  &
~\downarrow \nu \\
~\ \ \ \ \ \ \ \ \ \ \ \ \overset{\ast }{E} & ^{\underrightarrow{~\ \ h\circ
\overset{\ast }{\pi }~\ \ }} & N%
\end{array}%
\leqno(2.6)
\end{equation*}

Let $\Big(\overset{\overset{\ast }{\pi }^{\ast }\left( h^{\ast }F\right) }{%
\rho },Id_{E}\Big)$ be the\emph{\ }$\mathbf{B}^{\mathbf{v}}$-morphism of%
\emph{\ }$\left( \overset{\ast }{\pi }^{\ast }\left( h^{\ast }F\right) ,%
\overset{\ast }{\pi }^{\ast }\left( h^{\ast }\nu \right) ,\overset{\ast }{E}%
\right) $\ source and $\left( T\overset{\ast }{E},\tau _{\overset{\ast }{E}},%
\overset{\ast }{E}\right) $\ target, where%
\begin{equation*}
\begin{array}{rcl}
\overset{\ast }{\pi }^{\ast }\left( h^{\ast }F\right) & ^{\underrightarrow{%
\overset{\overset{\ast }{\pi }^{\ast }\left( h^{\ast }F\right) }{\rho }}} & T%
\overset{\ast }{E} \\
\displaystyle Z^{\alpha }T_{\alpha }\left( \overset{\ast }{u}_{x}\right) &
\longmapsto & \displaystyle\left( Z^{\alpha }\cdot \rho _{\alpha }^{i}\circ
h\circ \overset{\ast }{\pi }\right) \frac{\partial }{\partial x^{i}}\left(
\overset{\ast }{u}_{x}\right)%
\end{array}%
\leqno(2.7)
\end{equation*}

Using the operation
\begin{equation*}
\begin{array}{ccc}
\Gamma \left( \overset{\ast }{\pi }^{\ast }\left( h^{\ast }F\right) ,\overset%
{\ast }{\pi }^{\ast }\left( h^{\ast }\nu \right) ,\overset{\ast }{E}\right)
^{2} & ^{\underrightarrow{~\ \ \left[ ,\right] _{\overset{\ast }{\pi }^{\ast
}\left( h^{\ast }F\right) }~\ \ }} & \Gamma \left( \overset{\ast }{\pi }%
^{\ast }\left( h^{\ast }F\right) ,\overset{\ast }{\pi }^{\ast }\left(
h^{\ast }\nu \right) ,\overset{\ast }{E}\right)%
\end{array}%
\end{equation*}%
defined by%
\begin{equation*}
\begin{array}{ll}
\left[ T_{\alpha },T_{\beta }\right] _{\overset{\ast }{\pi }^{\ast }\left(
h^{\ast }F\right) } & =\left( L_{\alpha \beta }^{\gamma }\circ h\circ
\overset{\ast }{\pi }\right) T_{\gamma },\vspace*{1mm} \\
\left[ T_{\alpha },fT_{\beta }\right] _{\overset{\ast }{\pi }^{\ast }\left(
h^{\ast }F\right) } & \displaystyle=f\left( L_{\alpha \beta }^{\gamma }\circ
h\circ \overset{\ast }{\pi }\right) T_{\gamma }+\left( \rho _{\alpha
}^{i}\circ h\circ \overset{\ast }{\pi }\right) \frac{\partial f}{\partial
x^{i}}T_{\beta },\vspace*{1mm} \\
\left[ fT_{\alpha },T_{\beta }\right] _{\overset{\ast }{\pi }^{\ast }\left(
h^{\ast }F\right) } & =-\left[ T_{\beta },fT_{\alpha }\right] _{\overset{%
\ast }{\pi }^{\ast }\left( h^{\ast }F\right) },%
\end{array}%
\leqno(2.8)
\end{equation*}%
for any $f\in \mathcal{F}\left( \overset{\ast }{E}\right) ,$ it results that
\begin{equation*}
\begin{array}{c}
\left( \left( \overset{\ast }{\pi }^{\ast }\left( h^{\ast }F\right) ,\overset%
{\ast }{\pi }^{\ast }\left( h^{\ast }\nu \right) ,\overset{\ast }{E}\right) ,%
\left[ ,\right] _{\overset{\ast }{\pi }^{\ast }\left( h^{\ast }F\right)
},\left( \overset{\overset{\ast }{\pi }^{\ast }\left( h^{\ast }F\right) }{%
\rho },Id_{\overset{\ast }{E}}\right) \right)%
\end{array}%
\end{equation*}%
is a Lie algebroid.

\section{Natural and adapted basis}

In the following we consider the following diagram:
\begin{equation*}
\begin{array}{c}
\xymatrix{\overset{\ast }{E}\ar[d]_{\overset{\ast }{\pi }}&\left( F,\left[
,\right] _{F,h},\left( \rho ,\eta \right) \right)\ar[d]^\nu\\ M\ar[r]^h&N}%
\end{array}%
\leqno(3.1)
\end{equation*}%
where $\left( E,\pi ,M\right) \in \left\vert \mathbf{B}^{\mathbf{v}%
}\right\vert $ and $\left( \left( F,\nu ,N\right) ,\left[ ,\right]
_{F,h},\left( \rho ,\eta \right) \right) $ is a generalized Lie algebroid.

We take $\left( x^{i},p_{a}\right) $ as canonical local coordinates on $%
\left( \overset{\ast }{E},\overset{\ast }{\pi },M\right) ,$ where $i\in
\overline{1,m}$ and $a\in \overline{1,r}.$ Let
\begin{equation*}
\left( x^{i},p_{a}\right) \longrightarrow \left( x^{i%
{\acute{}}%
}\left( x^{i}\right) ,p_{a^{\prime }}\left( x^{i},p_{a}\right) \right)
\end{equation*}%
be a change of coordinates on $\left( \overset{\ast }{E},\overset{\ast }{\pi
},M\right) $. Then the coordinates $p_{a}$ change to $p_{a^{\prime }}$ by
the rule:
\begin{equation*}
\begin{array}{c}
p_{a^{\prime }}=M_{a^{\prime }}^{a}p_{a}.%
\end{array}%
\leqno(3.2)
\end{equation*}

Let
\begin{equation*}
\begin{array}[b]{c}
\left( \frac{\partial }{\partial x^{i}},\frac{\partial }{\partial p_{a}}%
\right) \overset{put}{=}\left( \overset{\ast }{\partial }_{i},\overset{\cdot
}{\partial }^{a}\right)%
\end{array}%
\leqno(3.3)
\end{equation*}%
be the natural base of the dual tangent Lie algebroid $\left( \left( T%
\overset{\ast }{E},\tau _{\overset{\ast }{E}},\overset{\ast }{E}\right) ,%
\left[ ,\right] _{T\overset{\ast }{E}},\left( Id_{T\overset{\ast }{E}},Id_{%
\overset{\ast }{E}}\right) \right) .$

For any sections%
\begin{equation*}
\begin{array}{c}
Z^{\alpha }T_{\alpha }\in \Gamma \left( \overset{\ast }{\pi }^{\ast }\left(
h^{\ast }F\right) ,\overset{\ast }{\pi }^{\ast }\left( h^{\ast }F\right) ,%
\overset{\ast }{E}\right)%
\end{array}%
\end{equation*}%
and%
\begin{equation*}
\begin{array}{c}
Y_{a}\displaystyle\overset{\cdot }{\partial }^{a}\in \Gamma \left( VT\overset%
{\ast }{E},\tau _{\overset{\ast }{E}},\overset{\ast }{E}\right)%
\end{array}%
\end{equation*}%
we obtain the section%
\begin{equation*}
\begin{array}{c}
Z^{\alpha }\overset{\ast }{\tilde{\partial}}_{\alpha }+Y_{a}\overset{\cdot }{%
\tilde{\partial}}^{a}=:Z^{\alpha }\left( T_{\alpha }\oplus \left( \rho
_{\alpha }^{i}\circ h\circ \overset{\ast }{\pi }\right) \overset{\ast }{%
\partial }_{i}\right) +Y_{a}\left( 0_{\overset{\ast }{\pi }^{\ast }\left(
h^{\ast }F\right) }\oplus \overset{\cdot }{\partial }^{a}\right) \vspace*{1mm%
} \\
=Z^{\alpha }T_{\alpha }\oplus \left( Z^{\alpha }\left( \rho _{\alpha
}^{i}\circ h\circ \overset{\ast }{\pi }\right) \overset{\ast }{\partial }%
_{i}+Y_{a}\overset{\cdot }{\partial }^{a}\right) \in \Gamma \left( \overset{%
\ast }{\pi }^{\ast }\left( h^{\ast }F\right) \oplus T\overset{\ast }{E},%
\overset{\oplus }{\pi },\overset{\ast }{E}\right) .%
\end{array}%
\end{equation*}

Since we have
\begin{equation*}
\begin{array}{c}
Z^{\alpha }\overset{\ast }{\tilde{\partial}}_{\alpha }+Y_{a}\overset{\cdot }{%
\tilde{\partial}}^{a}=0 \\
\Updownarrow \\
Z^{\alpha }T_{\alpha }=0~\wedge Z^{\alpha }\left( \rho _{\alpha }^{i}\circ
h\circ \overset{\ast }{\pi }\right) \overset{\ast }{\partial }_{i}+Y_{a}%
\overset{\cdot }{\partial }^{a}=0,%
\end{array}%
\end{equation*}%
it implies $Z^{\alpha }=0,~\alpha \in \overline{1,p}$ and $Y_{a}=0,~a\in
\overline{1,r}.$

Therefore, the sections $\displaystyle\overset{\ast }{\tilde{\partial}}%
_{1},...,\overset{\ast }{\tilde{\partial}}_{p},\overset{\cdot }{\tilde{%
\partial}}^{1},...,\overset{\cdot }{\tilde{\partial}}^{r}$ are linearly
independent.\smallskip

We consider the vector subbundle $\left( \left( \rho ,\eta \right) T\overset{%
\ast }{E},\left( \rho ,\eta \right) \tau _{\overset{\ast }{E}},\overset{\ast
}{E}\right) $ of the vector bundle\break $\left( \overset{\ast }{\pi }^{\ast
}\left( h^{\ast }F\right) \oplus T\overset{\ast }{E},\overset{\oplus }{\pi },%
\overset{\ast }{E}\right) ,$ for which the $\mathcal{F}\left( \overset{\ast }%
{E}\right) $-module of sections is the $\mathcal{F}\left( \overset{\ast }{E}%
\right) $-submodule of $\left( \Gamma \left( \overset{\ast }{\pi }^{\ast
}\left( h^{\ast }F\right) \oplus T\overset{\ast }{E},\overset{\oplus }{\pi },%
\overset{\ast }{E}\right) ,+,\cdot \right) ,$ generated by the set of
sections $\left( \overset{\ast }{\tilde{\partial}}_{\alpha },\overset{\cdot }%
{\tilde{\partial}}^{a}\right) $ which is called the \emph{natural }$\left(
\rho ,\eta \right) $\emph{-base.}

The matrix of coordinate transformation on $\left( \left( \rho ,\eta \right)
T\overset{\ast }{E},\left( \rho ,\eta \right) \tau _{\overset{\ast }{E}},%
\overset{\ast }{E}\right) $ at a change of fibred charts is
\begin{equation*}
\left\Vert
\begin{array}{cc}
\Lambda _{\alpha }^{\alpha
{\acute{}}%
}\circ h\circ \overset{\ast }{\pi } & 0\vspace*{1mm} \\
\left( \rho _{a}^{i}\circ h\circ \overset{\ast }{\pi }\right) \displaystyle%
\frac{\partial M_{b}^{a%
{\acute{}}%
}\circ \overset{\ast }{\pi }}{\partial x_{i}}y^{b} & M_{a}^{a%
{\acute{}}%
}\circ \overset{\ast }{\pi }%
\end{array}%
\right\Vert .\leqno(3.4)
\end{equation*}

We have the following

\textbf{Theorem 3.1 }\emph{Let} $\left( \overset{\ast }{\tilde{\rho}},Id_{%
\overset{\ast }{E}}\right) $\ \emph{be the} $\mathbf{B}^{\mathbf{v}}$\emph{%
-morphism of }$\left( \left( \rho ,\eta \right) T\overset{\ast }{E},\left(
\rho ,\eta \right) \tau _{\overset{\ast }{E}},\overset{\ast }{E}\right) $\
\emph{source and }$\left( T\overset{\ast }{E},\tau _{\overset{\ast }{E}},%
\overset{\ast }{E}\right) $\ \emph{target, where}
\begin{equation*}
\begin{array}{rcl}
\left( \rho ,\eta \right) T\overset{\ast }{E}\!\!\! & \!\!^{\underrightarrow{%
~\ \ \overset{\ast }{\tilde{\rho}}}~\ \ }\!\!\! & \!\!T\overset{\ast }{E}%
\vspace*{2mm} \\
\left( Z^{\alpha }\displaystyle\overset{\ast }{\tilde{\partial}}_{\alpha
}+Y_{a}\overset{\cdot }{\tilde{\partial}}^{a}\right) \!(\overset{\ast }{u}%
_{x})\!\!\!\! & \!\!\longmapsto \!\!\! & \!\!\left( \!Z^{\alpha }\!\left(
\rho _{\alpha }^{i}{\circ }h{\circ }\overset{\ast }{\pi }\!\right) \!\overset%
{\ast }{\partial }_{i}{+}Y_{a}\overset{\cdot }{\partial }^{a}\right) \!(%
\overset{\ast }{u}_{x})\!\!%
\end{array}%
\leqno(3.5)
\end{equation*}

\emph{Using the operation}
\begin{equation*}
\begin{array}{ccc}
\Gamma \left( \left( \rho ,\eta \right) T\overset{\ast }{E},\left( \rho
,\eta \right) \tau _{\overset{\ast }{E}},\overset{\ast }{E}\right) ^{2} & ^{%
\underrightarrow{~\ \ \left[ ,\right] _{\left( \rho ,\eta \right) T\overset{%
\ast }{E}}~\ \ }} & \Gamma \left( \left( \rho ,\eta \right) T\overset{\ast }{%
E},\left( \rho ,\eta \right) \tau _{\overset{\ast }{E}},\overset{\ast }{E}%
\right)%
\end{array}%
\end{equation*}%
\emph{defined by}%
\begin{equation*}
\begin{array}{l}
\left[ \left( Z_{1}^{\alpha }\displaystyle\overset{\ast }{\tilde{\partial}}%
_{\alpha }+Y_{a}^{1}\overset{\cdot }{\tilde{\partial}}^{a}\right) ,\left(
Z_{2}^{\beta }\displaystyle\overset{\ast }{\tilde{\partial}}_{\beta
}+Y_{b}^{2}\overset{\cdot }{\tilde{\partial}}^{b}\right) \right] _{\left(
\rho ,\eta \right) T\overset{\ast }{E}}\vspace*{1mm} \\
\displaystyle=\left[ Z_{1}^{\alpha }T_{a},Z_{2}^{\beta }T_{\beta }\right] _{%
\overset{\ast }{\pi }^{\ast }\left( h^{\ast }F\right) }\oplus \left[ \left(
\rho _{\alpha }^{i}\circ h\circ \overset{\ast }{\pi }\right) Z_{1}^{\alpha }%
\overset{\ast }{\partial }_{i}+Y_{a}^{1}\overset{\cdot }{\partial }%
^{a},\right. \vspace*{1mm} \\
\hfill \displaystyle\left. \left( \rho _{\beta }^{j}\circ h\circ \overset{%
\ast }{\pi }\right) Z_{2}^{\beta }\overset{\ast }{\partial }_{j}+Y_{b}^{2}%
\overset{\cdot }{\partial }^{b}\right] _{T\overset{\ast }{E}},%
\end{array}%
\leqno(3.6)
\end{equation*}%
\emph{for any} $\left( Z_{1}^{\alpha }\displaystyle\overset{\ast }{\tilde{%
\partial}}_{\alpha }+Y_{a}^{1}\overset{\cdot }{\tilde{\partial}}^{a}\right) $%
\emph{\ and }$\left( Z_{2}^{\beta }\displaystyle\overset{\ast }{\tilde{%
\partial}}_{\beta }+Y_{b}^{2}\overset{\cdot }{\tilde{\partial}}^{b}\right) ,$
\emph{we obtain that the couple }%
\begin{equation*}
\left( \left[ ,\right] _{\left( \rho ,\eta \right) T\overset{\ast }{E}%
},\left( \overset{\ast }{\tilde{\rho}},Id_{\overset{\ast }{E}}\right) \right)
\end{equation*}%
\emph{\ is a Lie algebroid structure for the vector bundle }$\left( \left(
\rho ,\eta \right) T\overset{\ast }{E},\left( \rho ,\eta \right) \tau _{%
\overset{\ast }{E}},\overset{\ast }{E}\right) .$

The Lie algebroid
\begin{equation*}
\begin{array}{c}
\left( \left( \left( \rho ,\eta \right) T\overset{\ast }{E},\left( \rho
,\eta \right) \tau _{\overset{\ast }{E}},\overset{\ast }{E}\right) ,\left[ ,%
\right] _{\left( \rho ,\eta \right) T\overset{\ast }{E}},\left( \overset{%
\ast }{\tilde{\rho}},Id_{\overset{\ast }{E}}\right) \right)%
\end{array}%
,
\end{equation*}%
is called the \emph{Lie algebroid generalized tangent bundle of dual vector
bundle }$\left( \overset{\ast }{E},\overset{\ast }{\pi },M\right) .$

We consider the $\mathbf{B}^{\mathbf{v}}$-morphism $\left( \left( \rho ,\eta
\right) \overset{\ast }{\pi }!,Id_{\overset{\ast }{E}}\right) $ given by the
commutative diagram%
\begin{equation*}
\begin{array}{c}
\xymatrix{\left( \rho ,\eta \right) T\overset{\ast }{E}\ar[r]^{( \rho ,\eta
) \overset{\ast }{\pi!} }\ar[d]_{(\rho,\eta)\tau_{\overset{\ast }{E}}}&
\overset{\ast }{\pi }^{\ast }\left( h^{\ast }F\right) \ar[d]^{pr_1} \\
\overset{\ast }{E}\ar[r]^{id_{\overset{\ast }{E}}}& \overset{\ast }{E}}%
\end{array}%
\leqno(3.7)
\end{equation*}

This is defined as:%
\begin{equation*}
\begin{array}[b]{c}
\left( \rho ,\eta \right) \overset{\ast }{\pi }!\left( \tilde{Z}^{\alpha }%
\overset{\ast }{\tilde{\partial}}_{\alpha }+Y_{a}\overset{\cdot }{\tilde{%
\partial}}^{a}\right) \left( \overset{\ast }{u}_{x}\right) =\left( \tilde{Z}%
^{\alpha }\tilde{T}_{\alpha }\right) \left( \overset{\ast }{u}_{x}\right) ,%
\end{array}%
\leqno(3.8)
\end{equation*}%
for any $\displaystyle\tilde{Z}^{\alpha }\overset{\ast }{\tilde{\partial}}%
_{\alpha }+Y_{a}\overset{\cdot }{\tilde{\partial}}^{a}\in \left( \left( \rho
,\eta \right) T\overset{\ast }{E},\left( \rho ,\eta \right) \tau _{\overset{%
\ast }{E}},\overset{\ast }{E}\right) .$\medskip

Using the $\mathbf{B}^{\mathbf{v}}$-morphisms $\left( 2.6\right) $ and $%
\left( 3.7\right) $ we obtain the \emph{tangent }$\left( \rho ,\eta \right) $%
\emph{-application }$\left( \left( \rho ,\eta \right) T\overset{\ast }{\pi }%
,h\circ \overset{\ast }{\pi }\right) $ of $\left( \left( \rho ,\eta \right) T%
\overset{\ast }{E},\left( \rho ,\eta \right) \tau _{\overset{\ast }{E}},%
\overset{\ast }{E}\right) $ source and $\left( F,\nu ,N\right) $ target.

\textbf{Definition 3.1} The kernel of the tangent $\left( \rho ,\eta \right)
$-application\ is written
\begin{equation*}
\left( V\left( \rho ,\eta \right) T\overset{\ast }{E},\left( \rho ,\eta
\right) \tau _{\overset{\ast }{E}},\overset{\ast }{E}\right)
\end{equation*}%
and it is called \emph{the vertical subbundle}.\bigskip

We remark that the set $\left\{ \displaystyle\overset{\cdot }{\tilde{\partial%
}}^{a},~a\in \overline{1,r}\right\} $ is a base of the $\mathcal{F}\left(
\overset{\ast }{E}\right) $-module
\begin{equation*}
\left( \Gamma \left( V\left( \rho ,\eta \right) T\overset{\ast }{E},\left(
\rho ,\eta \right) \tau _{\overset{\ast }{E}},\overset{\ast }{E}\right)
,+,\cdot \right) .
\end{equation*}

\textbf{Proposition 3.1} \emph{The short sequence of vector bundles }%
\begin{equation*}
\begin{array}{c}
\xymatrix{0\ar@{^(->}[r]^i\ar[d]&V(\rho,\eta)T\overset{\ast
}{E}\ar[d]\ar@{^(->}[r]^i&(\rho,\eta)T\overset{\ast
}{E}\ar[r]^{(\rho,\eta)\overset{\ast }{\pi }!}\ar[d]&\overset{\ast }{\pi
}^{\ast }\left( h^{\ast }F\right)\ar[r]\ar[d]&0\ar[d]\\ \overset{\ast
}{E}\ar[r]^{Id_{\overset{\ast }{E}}}&\overset{\ast
}{E}\ar[r]^{Id_{\overset{\ast }{E}}}&\overset{\ast
}{E}\ar[r]^{Id_{\overset{\ast }{E}}}&\overset{\ast }
{E}\ar[r]^{Id_{\overset{\ast }{E}}}&\overset{\ast }{E}}%
\end{array}%
\leqno(3.9)
\end{equation*}%
\emph{is exact.}

Let $\left( \rho ,\eta \right) \Gamma $ be a $\left( \rho ,\eta \right) $%
-connection for the vector bundle$\left( \overset{\ast }{E},\overset{\ast }{%
\pi },M\right) ,$ i. e. \textit{a }$\mathbf{Man}$-morphism of $\left( \rho
,\eta \right) T\overset{\ast }{E}$ source and $V\left( \rho ,\eta \right) T%
\overset{\ast }{E}$ target defined by%
\begin{equation*}
\begin{array}[b]{c}
\left( \rho ,\eta \right) \Gamma \left( \tilde{Z}^{\alpha }\overset{\ast }{%
\tilde{\partial}}_{\alpha }+Y_{b}\overset{\cdot }{\tilde{\partial}}%
^{b}\right) \left( \overset{\ast }{u}_{x}\right) =\left( Y_{b}-\left( \rho
,\eta \right) \Gamma _{b\alpha }\tilde{Z}^{\alpha }\right) \overset{\cdot }{%
\tilde{\partial}}^{b}\left( \overset{\ast }{u}_{x}\right) ,%
\end{array}%
\leqno(3.10)
\end{equation*}%
such that the $\mathbf{B}^{\mathbf{v}}$-morphism $\left( \left( \rho ,\eta
\right) \Gamma ,Id_{\overset{\ast }{E}}\right) $ is a split to the left in
the previous exact sequence. Its components satisfy the law of transformation%
\emph{\ }%
\begin{equation*}
\begin{array}[b]{c}
\left( \rho ,\eta \right) \Gamma _{b%
{\acute{}}%
\gamma
{\acute{}}%
}=M_{b%
{\acute{}}%
}^{b}{\circ }\overset{\ast }{\pi }\left[ -\left( \rho _{\gamma }^{i}\circ
h\circ \overset{\ast }{\pi }\right) \frac{\partial M_{b}^{a%
{\acute{}}%
}{\circ }\overset{\ast }{\pi }}{\partial x^{i}}p_{a%
{\acute{}}%
}+\left( \rho ,\eta \right) \Gamma _{b\gamma }\right] \left( \Lambda
_{\gamma
{\acute{}}%
}^{\gamma }\circ h\circ \overset{\ast }{\pi }\right) .%
\end{array}%
\leqno(3.11)
\end{equation*}

In the particular case of Lie algebroids, $\left( \eta ,h\right) =\left(
Id_{M},Id_{M}\right) ,$ we obtain\emph{\ }%
\begin{equation*}
\begin{array}[b]{c}
\rho \Gamma _{b%
{\acute{}}%
\gamma
{\acute{}}%
}=M_{b%
{\acute{}}%
}^{b}{\circ }\overset{\ast }{\pi }\left[ -\left( \rho _{\gamma }^{i}\circ
\overset{\ast }{\pi }\right) \frac{\partial M_{b}^{a%
{\acute{}}%
}{\circ }\overset{\ast }{\pi }}{\partial x^{i}}p_{a%
{\acute{}}%
}+\rho \Gamma _{b\gamma }\right] \left( \Lambda _{\gamma
{\acute{}}%
}^{\gamma }\circ \overset{\ast }{\pi }\right) .%
\end{array}%
\leqno(3.11)^{\prime }
\end{equation*}

In the classical case, $\left( \rho ,\eta ,h\right) =\left(
Id_{TE},Id_{M},Id_{M}\right) ,$ we obtain\emph{\ }%
\begin{equation*}
\begin{array}[b]{c}
\Gamma _{b%
{\acute{}}%
k%
{\acute{}}%
}=M_{b%
{\acute{}}%
}^{b}{\circ }\overset{\ast }{\pi }\left[ -\frac{\partial M_{b}^{a%
{\acute{}}%
}{\circ }\overset{\ast }{\pi }}{\partial x^{i}}p_{a%
{\acute{}}%
}+\Gamma _{bk}\right] \left( \frac{\partial x^{k}}{\partial x^{k^{\prime }}}%
\circ \overset{\ast }{\pi }\right) .%
\end{array}%
\leqno(3.11)^{\prime \prime }
\end{equation*}

The kernel of the $\mathbf{B}^{\mathbf{v}}$-morphism $\left( \left( \rho
,\eta \right) \Gamma ,Id_{\overset{\ast }{E}}\right) $\ is written $\left(
H\left( \rho ,\eta \right) T\overset{\ast }{E},\left( \rho ,\eta \right)
\tau _{\overset{\ast }{E}},\overset{\ast }{E}\right) $ and is called the
\emph{horizontal vector subbundle}.

We remark that the horizontal and the vertical vector subbundles are
interior differential systems (see $\left[ 2\right] $) of the Lie algebroid
generalized tangent bundle
\begin{equation*}
\begin{array}{c}
\left( \left( \left( \rho ,\eta \right) T\overset{\ast }{E},\left( \rho
,\eta \right) \tau _{\overset{\ast }{E}},\overset{\ast }{E}\right) ,\left[ ,%
\right] _{\left( \rho ,\eta \right) T\overset{\ast }{E}},\left( \overset{%
\ast }{\tilde{\rho}},Id_{\overset{\ast }{E}}\right) \right) .%
\end{array}%
\end{equation*}

We put the problem of finding a base for the $\mathcal{F}\left( \overset{%
\ast }{E}\right) $-module
\begin{equation*}
\left( \Gamma \left( H\left( \rho ,\eta \right) T\overset{\ast }{E},\left(
\rho ,\eta \right) \tau _{\overset{\ast }{E}},\overset{\ast }{E}\right)
,+,\cdot \right)
\end{equation*}%
of the type\textbf{\ }
\begin{equation*}
\begin{array}[t]{l}
\overset{\ast }{\tilde{\delta}}_{\alpha }=Z_{\alpha }^{\beta }\overset{\ast }%
{\tilde{\partial}}_{\beta }+Y_{a\alpha }\overset{\cdot }{\tilde{\partial}}%
^{a},\alpha \in \overline{1,r}%
\end{array}%
\end{equation*}%
which satisfies the following conditions:
\begin{equation*}
\begin{array}{rcl}
\displaystyle\Gamma \left( \left( \rho ,\eta \right) \overset{\ast }{\pi }%
!,Id_{\overset{\ast }{E}}\right) \left( \overset{\ast }{\tilde{\delta}}%
_{\alpha }\right) & = & T_{\alpha }\vspace*{2mm}, \\
\displaystyle\Gamma \left( \left( \rho ,\eta \right) \Gamma ,Id_{\overset{%
\ast }{E}}\right) \left( \overset{\ast }{\tilde{\delta}}_{\alpha }\right) & =
& 0.%
\end{array}%
\leqno(3.12)
\end{equation*}

Then we obtain the sections
\begin{equation*}
\begin{array}[t]{l}
\displaystyle\overset{\ast }{\tilde{\delta}}_{\alpha }=\overset{\ast }{%
\tilde{\partial}}_{\alpha }+\left( \rho ,\eta \right) \Gamma _{b\alpha }%
\overset{\cdot }{\tilde{\partial}}^{b}=T_{\alpha }\oplus \left( \left( \rho
_{\alpha }^{i}\circ h\circ \overset{\ast }{\pi }\right) \overset{\ast }{%
\partial }_{i}+\left( \rho ,\eta \right) \Gamma _{b\alpha }\overset{\cdot }{%
\partial }^{b}\right) .%
\end{array}%
\leqno(3.13)
\end{equation*}%
such that their law of change is a tensorial law under a change of vector
fiber charts.

The base $\left( \overset{\ast }{\tilde{\delta}}_{\alpha },\overset{\cdot }{%
\tilde{\partial}}^{a}\right) $ will be called the \emph{adapted }$\left(
\rho ,\eta \right) $\emph{-base.}

\textit{Remark 3.2 }The following equality holds good%
\begin{equation*}
\begin{array}{l}
\Gamma \left( \overset{\ast }{\tilde{\rho}},Id_{\overset{\ast }{E}}\right)
\left( \overset{\ast }{\tilde{\delta}}_{\alpha }\right) =\left( \rho
_{\alpha }^{i}\circ h\circ \overset{\ast }{\pi }\right) \overset{\ast }{%
\partial }_{i}+\left( \rho ,\eta \right) \Gamma _{b\alpha }\dot{\partial}%
^{b}.%
\end{array}%
\leqno(3.14)
\end{equation*}

Moreover, if $\left( \rho ,\eta \right) \Gamma $ is the $\left( \rho ,\eta
\right) $-connection associated to a connection $\Gamma $ (see $\left[ 1%
\right] $), then we obtain
\begin{equation*}
\begin{array}{l}
\Gamma \left( \overset{\ast }{\tilde{\rho}},Id_{\overset{\ast }{E}}\right)
\left( \overset{\ast }{\tilde{\delta}}_{\alpha }\right) =\left( \rho
_{\alpha }^{i}\circ h\circ \overset{\ast }{\pi }\right) \overset{\ast }{%
\delta }_{i},%
\end{array}%
\leqno(3.15)
\end{equation*}%
where $\left( \overset{\ast }{\delta }_{i},\dot{\partial}^{a}\right) $ is
the adapted base for the $\mathcal{F}\left( \overset{\ast }{E}\right) $%
-module $\left( \Gamma \left( T\overset{\ast }{E},\tau _{\overset{\ast }{E}},%
\overset{\ast }{E}\right) ,+,\cdot \right) .$

Let $\left( d\tilde{z}^{\alpha },d\tilde{p}_{a}\right) $ be the natural dual
$\left( \rho ,\eta \right) $-base of natural $\left( \rho ,\eta \right) $%
-base $\left( \displaystyle\overset{\ast }{\partial }_{\alpha },\displaystyle%
\overset{\cdot }{\tilde{\partial}}^{a}\right) .$

This is determined by the equations
\begin{equation*}
\begin{array}{c}
\left\{
\begin{array}{cc}
\displaystyle\left\langle d\tilde{z}^{\alpha },\overset{\ast }{\tilde{%
\partial}}_{\beta }\right\rangle =\delta _{\beta }^{\alpha }, & \displaystyle%
\left\langle d\tilde{z}^{\alpha },\overset{\cdot }{\tilde{\partial}}%
^{b}\right\rangle =0,\vspace*{2mm} \\
\displaystyle\left\langle d\tilde{p}_{a},\overset{\ast }{\tilde{\partial}}%
_{\beta }\right\rangle =0, & \displaystyle\left\langle d\tilde{p}_{a},%
\overset{\cdot }{\tilde{\partial}}^{b}\right\rangle =\delta _{a}^{b}.%
\end{array}%
\right.%
\end{array}%
\end{equation*}

We consider the problem of finding a base for the $\mathcal{F}\left( \overset%
{\ast }{E}\right) $-module
\begin{equation*}
\left( \Gamma \left( \left( V\left( \rho ,\eta \right) T\overset{\ast }{E}%
\right) ^{\ast },\left( \left( \rho ,\eta \right) \tau _{\overset{\ast }{E}%
}\right) ^{\ast },\overset{\ast }{E}\right) ,+,\cdot \right)
\end{equation*}%
of the type
\begin{equation*}
\begin{array}{c}
\delta \tilde{p}_{a}=\theta _{a\alpha }d\tilde{z}^{\alpha }+\omega _{a}^{b}d%
\tilde{p}_{b},~a\in \overline{1,r}%
\end{array}%
\end{equation*}%
which satisfies the following conditions:
\begin{equation*}
\begin{array}{c}
\left\langle \delta \tilde{p}_{a},\overset{\cdot }{\tilde{\partial}}%
^{b}\right\rangle =\delta _{a}^{b}\wedge \left\langle \delta \tilde{p}_{a},%
\overset{\ast }{\tilde{\delta}}_{\alpha }\right\rangle =0,%
\end{array}%
\leqno(3.16)
\end{equation*}

We obtain the sections
\begin{equation*}
\begin{array}{l}
\delta \tilde{p}_{a}=-\left( \rho ,\eta \right) \Gamma _{a\alpha }d\tilde{z}%
^{\alpha }+d\tilde{p}_{a},a\in \overline{1,r}.%
\end{array}%
\leqno(3.17)
\end{equation*}%
such that their changing rule is tensorial under a change of vector fiber
charts. The base $\left( d\tilde{z}^{\alpha },\delta \tilde{p}_{a}\right) $
will be called the \emph{adapted dual }$\left( \rho ,\eta \right) $\emph{%
-base.}

\section{Tensor $d$-fields. Distinguished linear $\left( \protect\rho ,%
\protect\eta \right) $-connections}

We consider the following diagram:
\begin{equation*}
\begin{array}{c}
\xymatrix{\overset{\ast }{E}\ar[d]_{\overset{\ast }{\pi }}&\left( F,\left[
,\right] _{F,h},\left( \rho ,\eta \right) \right)\ar[d]^\nu\\ M\ar[r]^h&N}%
\end{array}%
\end{equation*}%
where $\left( E,\pi ,M\right) \in \left\vert \mathbf{B}^{\mathbf{v}%
}\right\vert $ and $\left( \left( F,\nu ,N\right) ,\left[ ,\right]
_{F,h},\left( \rho ,\eta \right) \right) $ is a generalized Lie algebroid.

Let $\left( \rho ,\eta \right) \Gamma $ be a $\left( \rho ,\eta \right) $%
-connection $\ $for the vector bundle $\left( \overset{\ast }{E},\overset{%
\ast }{\pi },M\right) .$

Let
\begin{equation*}
\left( \mathcal{T}~_{q,s}^{p,r}\left( \left( \rho ,\eta \right) T\overset{%
\ast }{E},\left( \rho ,\eta \right) \tau _{\overset{\ast }{E}},\overset{\ast
}{E}\right) ,+,\cdot \right)
\end{equation*}%
be the $\mathcal{F}\left( \overset{\ast }{E}\right) $-module of tensor
fields by $\left( _{q,s}^{p,r}\right) $-type from the generalized tangent
bundle
\begin{equation*}
\left( H\left( \rho ,\eta \right) T\overset{\ast }{E}\oplus V\left( \rho
,\eta \right) T\overset{\ast }{E},\left( \rho ,\eta \right) \tau _{\overset{%
\ast }{E}},\overset{\ast }{E}\right) .
\end{equation*}

An arbitrarily tensor field $T$ is written by the form:
\begin{equation*}
T=T_{\beta _{1}...\beta _{q}b_{1}...b_{s}}^{\alpha _{1}...\alpha
_{p}a_{1}...a_{r}}\overset{\ast }{\tilde{\delta}}_{\alpha _{1}}\otimes
...\otimes \overset{\ast }{\tilde{\delta}}_{\alpha _{p}}\otimes d\tilde{z}%
^{\beta _{1}}\otimes ...\otimes d\tilde{z}^{\beta _{q}}\otimes \overset{%
\cdot }{\tilde{\partial}}^{b_{1}}\otimes ...\otimes \overset{\cdot }{\tilde{%
\partial}}^{b_{s}}\otimes \delta \tilde{p}_{a_{1}}\otimes ...\otimes \delta
\tilde{p}_{a_{r}}.
\end{equation*}

Let
\begin{equation*}
\left( ~\mathcal{T}\left( \left( \rho ,\eta \right) T\overset{\ast }{E}%
,\left( \rho ,\eta \right) \tau _{\overset{\ast }{E}},\overset{\ast }{E}%
\right) ,+,\cdot ,\otimes \right)
\end{equation*}%
be the tensor fields algebra of generalized tangent bundle $\left( \left(
\rho ,\eta \right) T\overset{\ast }{E},\left( \rho ,\eta \right) \tau _{%
\overset{\ast }{E}},\overset{\ast }{E}\right) $.

If $T_{1}{\in }\mathcal{T}_{q_{1},s_{1}}^{p_{1},r_{1}}\left( \left( \rho
,\eta \right) T\overset{\ast }{E},\left( \rho ,\eta \right) \tau _{\overset{%
\ast }{E}},\overset{\ast }{E}\right) $ and $T_{2}{\in }\mathcal{T}%
_{q_{2},s_{2}}^{p_{2},r_{2}}\left( \left( \rho ,\eta \right) T\overset{\ast }%
{E},\left( \rho ,\eta \right) \tau _{\overset{\ast }{E}},\overset{\ast }{E}%
\right) $, then the components of product tensor field $T_{1}\otimes T_{2}$
are the products of local components of $T_{1}$ and $T_{2}.$

Therefore, we obtain $T_{1}\otimes T_{2}\in \mathcal{T}%
_{q_{1}+q_{2},s_{1}+s_{2}}^{p_{1}+p_{2},r_{1}+r_{2}}\left( \left( \rho ,\eta
\right) T\overset{\ast }{E},\left( \rho ,\eta \right) \tau _{\overset{\ast }{%
E}},\overset{\ast }{E}\right) .$

Let $\mathcal{DT}\left( \left( \rho ,\eta \right) T\overset{\ast }{E},\left(
\rho ,\eta \right) \tau _{\overset{\ast }{E}},\overset{\ast }{E}\right) $ be
the family of tensor fields
\begin{equation*}
T\in \mathcal{T}\left( \left( \rho ,\eta \right) T\overset{\ast }{E},\left(
\rho ,\eta \right) \tau _{\overset{\ast }{E}},\overset{\ast }{E}\right)
\end{equation*}%
for which there exists%
\begin{equation*}
T_{1}{\in }\mathcal{T}_{q,0}^{p,0}\left( \left( \rho ,\eta \right) T\overset{%
\ast }{E},\left( \rho ,\eta \right) \tau _{\overset{\ast }{E}},\overset{\ast
}{E}\right)
\end{equation*}%
and%
\begin{equation*}
T_{2}{\in }\mathcal{T}_{0,s}^{0,r}\left( \left( \rho ,\eta \right) T\overset{%
\ast }{E},\left( \rho ,\eta \right) \tau _{\overset{\ast }{E}},\overset{\ast
}{E}\right)
\end{equation*}
such that $T=T_{1}+T_{2}.$

The $\mathcal{F}\left( \overset{\ast }{E}\right) $-module $\left( \mathcal{DT%
}\left( \left( \rho ,\eta \right) T\overset{\ast }{E},\left( \rho ,\eta
\right) \tau _{\overset{\ast }{E}},\overset{\ast }{E}\right) ,+,\cdot
\right) $ will be called the \emph{module of distinguished tensor fields} or
the \emph{module of tensor }$d$-\emph{fields.}

\emph{\ Remark 4.1 }The elements of
\begin{equation*}
\Gamma \left( \left( \rho ,\eta \right) T\overset{\ast }{E},\left( \rho
,\eta \right) \tau _{\overset{\ast }{E}},\overset{\ast }{E}\right)
\end{equation*}%
respectively
\begin{equation*}
\Gamma (((\rho ,\eta )T\overset{\ast }{E})^{\ast },\break ((\rho ,\eta )\tau
_{\overset{\ast }{E}})^{\ast },\overset{\ast }{E})
\end{equation*}%
are tensor $d$-fields.

\textbf{Definition 4.1 }Let $\left( \rho ,\eta \right) \Gamma $ be a $\left(
\rho ,\eta \right) $-connection for the vector bundle $\left( \overset{\ast }%
{E},\overset{\ast }{\pi },M\right) $ and let
\begin{equation*}
\begin{array}{l}
\left( X,T\right) ^{\underrightarrow{\left( \rho ,\eta \right) \overset{\ast
}{D}}\,}\vspace*{1mm}\left( \rho ,\eta \right) \overset{\ast }{D}_{X}T%
\end{array}%
\leqno(4.1)
\end{equation*}%
be a covariant $\left( \rho ,\eta \right) $-derivative for the tensor
algebra of generalized tangent bundle
\begin{equation*}
\left( \left( \rho ,\eta \right) T\overset{\ast }{E},\left( \rho ,\eta
\right) \tau _{\overset{\ast }{E}},\overset{\ast }{E}\right)
\end{equation*}%
which \ preserves \ the \ horizontal and vertical \emph{IDS} by parallelism.
(see $\left[ 2\right] $)

If $\left( U,\overset{\ast }{s}_{U}\right) $ is a vector local $\left(
m+r\right) $-chart for $\left( \overset{\ast }{E},\overset{\ast }{\pi }%
,M\right) ,$ then the real local functions
\begin{equation*}
\left( \left( \rho ,\eta \right) \overset{\ast }{H}_{\beta \gamma }^{\alpha
},\left( \rho ,\eta \right) \overset{\ast }{H}_{b\gamma }^{a},\left( \rho
,\eta \right) \overset{\ast }{V}_{\beta }^{\alpha c},\left( \rho ,\eta
\right) \overset{\ast }{V}_{a}^{bc}\right)
\end{equation*}%
defined on $\overset{\ast }{\pi }^{-1}\left( U\right) $ and determined by
the following equalities:
\begin{equation*}
\begin{array}{ll}
\left( \rho ,\eta \right) \overset{\ast }{D}_{\overset{\ast }{\tilde{\delta}}%
_{\gamma }}\overset{\ast }{\tilde{\delta}}_{\beta }=\left( \rho ,\eta
\right) \overset{\ast }{H}_{\beta \gamma }^{\alpha }\overset{\ast }{\tilde{%
\delta}}_{\alpha }, & \left( \rho ,\eta \right) \overset{\ast }{D}_{\overset{%
\ast }{\tilde{\delta}}_{\gamma }}\overset{\cdot }{\tilde{\partial}}%
^{a}=\left( \rho ,\eta \right) \overset{\ast }{H}_{b\gamma }^{a}\overset{%
\cdot }{\tilde{\partial}}^{b} \\
\left( \rho ,\eta \right) \overset{\ast }{D}_{\overset{\cdot }{\tilde{%
\partial}}^{c}}\overset{\ast }{\tilde{\delta}}_{\beta }=\left( \rho ,\eta
\right) \overset{\ast }{V}_{\beta }^{\alpha c}\overset{\ast }{\tilde{\delta}}%
_{\alpha }, & \left( \rho ,\eta \right) \overset{\ast }{D}_{\overset{\cdot }{%
\tilde{\partial}}^{c}}\overset{\cdot }{\tilde{\partial}}^{b}=\left( \rho
,\eta \right) \overset{\ast }{V}_{a}^{bc}\overset{\cdot }{\tilde{\partial}}%
^{a}%
\end{array}%
\leqno(4.2)
\end{equation*}%
are the components of a linear $\left( \rho ,\eta \right) $-connection
\begin{equation*}
\left( \left( \rho ,\eta \right) \overset{\ast }{H},\left( \rho ,\eta
\right) \overset{\ast }{V}\right)
\end{equation*}%
for the generalized tangent bundle $\left( \left( \rho ,\eta \right) T%
\overset{\ast }{E},\left( \rho ,\eta \right) \tau _{\overset{\ast }{E}},%
\overset{\ast }{E}\right) $ which will be called the \emph{distinguished
linear }$\left( \rho ,\eta \right) $\emph{-connection.}

If $h=Id_{M},$ then the distinguished linear $\left( Id_{TM},Id_{M}\right) $%
-connection will be called the \emph{distinguished linear connection.}

The components of a distinguished linear connection $\left( \overset{\ast }{H%
},\overset{\ast }{V}\right) $ will be denoted
\begin{equation*}
\left( \overset{\ast }{H}_{jk}^{i},\overset{\ast }{H}_{bk}^{a},\overset{\ast
}{V}_{j}^{ic},\overset{\ast }{V}_{a}^{bc}\right) .
\end{equation*}

\textbf{Theorem 4.1 }\emph{If }$\left( \left( \rho ,\eta \right) \overset{%
\ast }{H},\left( \rho ,\eta \right) \overset{\ast }{V}\right) $ \emph{is a
distinguished linear} $(\rho ,\eta )$-\emph{connection for the generalized
tangent bundle }$\left( \left( \rho ,\eta \right) T\overset{\ast }{E},\left(
\rho ,\eta \right) \tau _{\overset{\ast }{E}},\overset{\ast }{E}\right) $%
\emph{, then its components satisfy the change relations: }

\begin{equation*}
\begin{array}{ll}
\left( \rho ,\eta \right) \overset{\ast }{H}_{\beta
{\acute{}}%
\gamma
{\acute{}}%
}^{\alpha
{\acute{}}%
}\!\! & =\Lambda _{\alpha }^{\alpha
{\acute{}}%
}\circ h\circ \overset{\ast }{\pi }\left[ \Gamma \left( \overset{\ast }{%
\tilde{\rho}},Id_{\overset{\ast }{E}}\right) \left( \overset{\ast }{\tilde{%
\delta}}_{\gamma }\right) \left( \Lambda _{\beta
{\acute{}}%
}^{\alpha }\circ h\circ \overset{\ast }{\pi }\right) +\right. \vspace*{1mm}
\\
& +\left. \left( \rho ,\eta \right) \overset{\ast }{H}_{\beta \gamma
}^{\alpha }\cdot \Lambda _{\beta
{\acute{}}%
}^{\beta }\circ h\circ \overset{\ast }{\pi }\right] \cdot \Lambda _{\gamma
{\acute{}}%
}^{\gamma }\circ h\circ \overset{\ast }{\pi },\vspace*{2mm} \\
\left( \rho ,\eta \right) \overset{\ast }{H}_{b%
{\acute{}}%
\gamma
{\acute{}}%
}^{a%
{\acute{}}%
}\!\! & =M_{a}^{a%
{\acute{}}%
}\circ \overset{\ast }{\pi }\left[ \Gamma \left( \overset{\ast }{\tilde{\rho}%
},Id_{\overset{\ast }{E}}\right) \left( \overset{\ast }{\tilde{\delta}}%
_{\gamma }\right) \left( M_{b%
{\acute{}}%
}^{a}\circ \overset{\ast }{\pi }\right) +\right. \vspace*{1mm} \\
& \left. +\left( \rho ,\eta \right) \overset{\ast }{H}_{b\gamma }^{a}\cdot
M_{b%
{\acute{}}%
}^{b}\circ \overset{\ast }{\pi }\right] \cdot \Lambda _{\gamma
{\acute{}}%
}^{\gamma }\circ h\circ \overset{\ast }{\pi },\vspace*{2mm} \\
\left( \rho ,\eta \right) \overset{\ast }{V}_{\beta
{\acute{}}%
}^{\alpha
{\acute{}}%
c%
{\acute{}}%
}\!\! & =\Lambda _{\alpha }^{\alpha
{\acute{}}%
}\circ h\circ \overset{\ast }{\pi }\cdot \left( \rho ,\eta \right) \overset{%
\ast }{V}_{\beta }^{\alpha c}\cdot \Lambda _{\beta
{\acute{}}%
}^{\beta }\circ h\circ \overset{\ast }{\pi }\cdot M_{c}^{c%
{\acute{}}%
}\circ \overset{\ast }{\pi },\vspace*{2mm} \\
\left( \rho ,\eta \right) \overset{\ast }{V}_{b%
{\acute{}}%
}^{a%
{\acute{}}%
c%
{\acute{}}%
}\!\! & =M_{a}^{a%
{\acute{}}%
}\circ \overset{\ast }{\pi }\cdot \left( \rho ,\eta \right) \overset{\ast }{V%
}_{b}^{ac}\cdot M_{b%
{\acute{}}%
}^{b}\circ \overset{\ast }{\pi }\cdot M_{c}^{c%
{\acute{}}%
}\circ \overset{\ast }{\pi }.%
\end{array}%
\leqno(4.3)
\end{equation*}

\textbf{Corollary 4.1 }\emph{In the particular case of Lie algebroids, }$%
\left( \eta ,h\right) =\left( Id_{M},Id_{M}\right) ,$\emph{\ we obtain}%
\begin{equation*}
\begin{array}{ll}
\rho \overset{\ast }{H}_{\beta
{\acute{}}%
\gamma
{\acute{}}%
}^{\alpha
{\acute{}}%
}\!\! & =\Lambda _{\alpha }^{\alpha
{\acute{}}%
}\circ \overset{\ast }{\pi }\left[ \Gamma \left( \overset{\ast }{\tilde{\rho}%
},Id_{\overset{\ast }{E}}\right) \left( \overset{\ast }{\tilde{\delta}}%
_{\gamma }\right) \left( \Lambda _{\beta
{\acute{}}%
}^{\alpha }\circ \overset{\ast }{\pi }\right) +\rho \overset{\ast }{H}%
_{\beta \gamma }^{\alpha }\cdot \Lambda _{\beta
{\acute{}}%
}^{\beta }\circ \overset{\ast }{\pi }\right] \cdot \Lambda _{\gamma
{\acute{}}%
}^{\gamma }\circ \overset{\ast }{\pi },\vspace*{2mm} \\
\rho \overset{\ast }{H}_{b%
{\acute{}}%
\gamma
{\acute{}}%
}^{a%
{\acute{}}%
}\!\! & =M_{a}^{a%
{\acute{}}%
}\circ \overset{\ast }{\pi }\left[ \Gamma \left( \overset{\ast }{\tilde{\rho}%
},Id_{\overset{\ast }{E}}\right) \left( \overset{\ast }{\tilde{\delta}}%
_{\gamma }\right) \left( M_{b%
{\acute{}}%
}^{a}\circ \overset{\ast }{\pi }\right) +\rho \overset{\ast }{H}_{b\gamma
}^{a}\cdot M_{b%
{\acute{}}%
}^{b}\circ \overset{\ast }{\pi }\right] \cdot \Lambda _{\gamma
{\acute{}}%
}^{\gamma }\circ \overset{\ast }{\pi },\vspace*{2mm} \\
\rho \overset{\ast }{V}_{\beta
{\acute{}}%
}^{\alpha
{\acute{}}%
c%
{\acute{}}%
}\!\! & =\Lambda _{\alpha }^{\alpha
{\acute{}}%
}\circ \overset{\ast }{\pi }\cdot \rho \overset{\ast }{V}_{\beta }^{\alpha
c}\cdot \Lambda _{\beta
{\acute{}}%
}^{\beta }\circ \overset{\ast }{\pi }\cdot M_{c}^{c%
{\acute{}}%
}\circ \overset{\ast }{\pi },\vspace*{2mm} \\
\rho \overset{\ast }{V}_{b%
{\acute{}}%
}^{a%
{\acute{}}%
c%
{\acute{}}%
}\!\! & =M_{a}^{a%
{\acute{}}%
}\circ \overset{\ast }{\pi }\cdot \rho \overset{\ast }{V}_{b}^{ac}\cdot M_{b%
{\acute{}}%
}^{b}\circ \overset{\ast }{\pi }\cdot M_{c}^{c%
{\acute{}}%
}\circ \overset{\ast }{\pi }.%
\end{array}%
\leqno(4.3)^{\prime }
\end{equation*}

\emph{In the classical case, }$\left( \rho ,\eta ,h\right) =\left(
Id_{TE},Id_{M},Id_{M}\right) ,$\emph{\ we obtain that the components of a
distinguished linear connection $\left( H,V\right) $ verify the change
relations:}
\begin{equation*}
\begin{array}{cl}
\overset{\ast }{H}_{j%
{\acute{}}%
k%
{\acute{}}%
}^{i%
{\acute{}}%
} & =\frac{\partial x^{i%
{\acute{}}%
}}{\partial x^{i}}\circ \overset{\ast }{\pi }\cdot \left[ \frac{\overset{%
\ast }{\delta }}{\delta x^{k}}\left( \frac{\partial x^{i}}{\partial x^{j%
{\acute{}}%
}}\circ \overset{\ast }{\pi }\right) +\overset{\ast }{H}_{jk}^{i}\cdot \frac{%
\partial x^{j}}{\partial x^{j%
{\acute{}}%
}}\circ \overset{\ast }{\pi }\right] \cdot \frac{\partial x^{k}}{\partial
x^{k%
{\acute{}}%
}}\circ \overset{\ast }{\pi },\vspace*{2mm} \\
\overset{\ast }{H}_{b%
{\acute{}}%
k%
{\acute{}}%
}^{a%
{\acute{}}%
} & =M_{a}^{a%
{\acute{}}%
}\circ \overset{\ast }{\pi }\cdot \left[ \frac{\overset{\ast }{\delta }}{%
\delta x^{k}}\left( M_{b%
{\acute{}}%
}^{a}\circ \overset{\ast }{\pi }\right) +\overset{\ast }{H}_{bk}^{a}\cdot
M_{b%
{\acute{}}%
}^{b}\circ \overset{\ast }{\pi }\right] \cdot \frac{\partial x^{k}}{\partial
x^{k%
{\acute{}}%
}}\circ \overset{\ast }{\pi },\vspace*{2mm} \\
\overset{\ast }{V}_{j%
{\acute{}}%
}^{i%
{\acute{}}%
c%
{\acute{}}%
} & =\frac{\partial x^{i%
{\acute{}}%
}}{\partial x^{i}}\circ \overset{\ast }{\pi }\cdot \overset{\ast }{V}%
_{j}^{ic}\frac{\partial x^{j}}{\partial x^{j%
{\acute{}}%
}}\circ \overset{\ast }{\pi }\cdot M_{c}^{c^{\prime }}\circ \overset{\ast }{%
\pi },\vspace*{3mm} \\
\overset{\ast }{V}_{b%
{\acute{}}%
}^{a%
{\acute{}}%
c%
{\acute{}}%
} & =M_{a}^{a%
{\acute{}}%
}\circ \overset{\ast }{\pi }\cdot \overset{\ast }{V}_{b}^{ac}\cdot M_{b%
{\acute{}}%
}^{b}\circ \overset{\ast }{\pi }\cdot M_{c}^{c%
{\acute{}}%
}\circ \overset{\ast }{\pi }.%
\end{array}%
\leqno(4.3)^{\prime \prime }
\end{equation*}

\textbf{Example 4.1 }If $\left( \overset{\ast }{E},\overset{\ast }{\pi }%
,M\right) $ is endowed with the $\left( \rho ,\eta \right) $-connection $%
\left( \rho ,\eta \right) \Gamma $, then the local real functions
\begin{equation*}
\begin{array}[b]{c}
\left( \frac{\partial \left( \rho ,\eta \right) \Gamma _{b\gamma }}{\partial
p_{a}},\frac{\partial \left( \rho ,\eta \right) \Gamma _{b\gamma }}{\partial
p_{a}},0,0\right)%
\end{array}%
\leqno(4.4)
\end{equation*}%
are the components of a distinguished linear $\left( \rho ,\eta \right) $%
\textit{-}connection for the generalized tangent bundle
\begin{equation*}
\left( \left( \rho ,\eta \right) T\overset{\ast }{E},\left( \rho ,\eta
\right) \tau _{\overset{\ast }{E}},\overset{\ast }{E}\right) ,
\end{equation*}%
which will by called the \emph{Berwald linear }$\left( \rho ,\eta \right) $%
\emph{-connection.}

\textbf{Theorem 4.2} \emph{If the generalized tangent bundle} $\!\left(
\left( \rho ,\eta \right) T\overset{\ast }{E},\left( \rho ,\eta \right) \tau
_{\overset{\ast }{E}},\overset{\ast }{E}\right) $ \emph{is endowed with a
distinguished linear} $\!(\rho ,\!\eta )$\emph{-connection} $((\rho ,\eta )%
\overset{\ast }{H},(\rho ,\eta )\overset{\ast }{V})$, \emph{then, for any}
\begin{equation*}
X=Z^{\gamma }\overset{\ast }{\tilde{\delta}}_{\gamma }+Y_{a}\overset{\cdot }{%
\tilde{\partial}}^{a}\in \Gamma \left( \left( \rho ,\eta \right) T\overset{%
\ast }{E},\left( \rho ,\eta \right) \tau _{\overset{\ast }{E}},\overset{\ast
}{E}\right)
\end{equation*}%
\emph{and for any}
\begin{equation*}
T\in \mathcal{T}_{qs}^{pr}\!\left( \left( \rho ,\eta \right) T\overset{\ast }%
{E},\left( \rho ,\eta \right) \tau _{\overset{\ast }{E}},\overset{\ast }{E}%
\right) ,
\end{equation*}%
\emph{we obtain the formula:}
\begin{equation*}
\begin{array}{l}
\left( \rho ,\eta \right) D_{X}\left( T_{\beta _{1}...\beta
_{q}b_{1}...b_{s}}^{\alpha _{1}...\alpha _{p}a_{1}...a_{r}}\overset{\ast }{%
\tilde{\delta}}_{\alpha _{1}}\otimes ...\otimes \overset{\ast }{\tilde{\delta%
}}_{\alpha _{p}}\otimes d\tilde{z}^{\beta _{1}}\otimes ...\otimes \right.
\vspace*{1mm} \\
\hspace*{9mm}\left. \otimes d\tilde{z}^{\beta _{q}}\otimes \overset{\cdot }{%
\tilde{\partial}}^{b_{1}}\otimes ...\otimes \overset{\cdot }{\tilde{\partial}%
}^{b_{s}}\otimes \delta \tilde{p}_{a_{1}}\otimes ...\otimes \delta \tilde{p}%
_{a_{r}}\right) =\vspace*{1mm} \\
\hspace*{9mm}=Z^{\gamma }T_{\beta _{1}...\beta _{q}b_{1}...b_{s}\mid \gamma
}^{\alpha _{1}...\alpha _{p}a_{1}...a_{r}}\overset{\ast }{\tilde{\delta}}%
_{\alpha _{1}}\otimes ...\otimes \overset{\ast }{\tilde{\delta}}_{\alpha
_{p}}\otimes d\tilde{z}^{\beta _{1}}\otimes ...\otimes d\tilde{z}^{\beta
_{q}}\otimes \overset{\cdot }{\tilde{\partial}}^{b_{1}}\otimes ...\otimes
\vspace*{1mm} \\
\hspace*{9mm}\otimes \overset{\cdot }{\tilde{\partial}}^{b_{s}}\otimes
\delta \tilde{p}_{a_{1}}\otimes ...\otimes \delta \tilde{p}%
_{a_{r}}+Y_{c}T_{\beta _{1}...\beta _{q}b_{1}...b_{s}}^{\alpha _{1}...\alpha
_{p}a_{1}...a_{r}}\mid ^{c}\overset{\ast }{\tilde{\delta}}_{\alpha
_{1}}\otimes ...\otimes \vspace*{1mm} \\
\hspace*{9mm}\otimes \overset{\ast }{\tilde{\delta}}_{\alpha _{p}}\otimes d%
\tilde{z}^{\beta _{1}}\otimes ...\otimes d\tilde{z}^{\beta _{q}}\otimes
\overset{\cdot }{\tilde{\partial}}^{b_{1}}\otimes ...\otimes \overset{\cdot }%
{\tilde{\partial}}^{b_{s}}\otimes \delta \tilde{p}_{a_{1}}\otimes ...\otimes
\delta \tilde{p}_{a_{r}},%
\end{array}%
\leqno(4.5)
\end{equation*}%
\emph{where\ }%
\begin{equation*}
\begin{array}{l}
T_{\beta _{1}...\beta _{q}b_{1}...b_{s}\mid \gamma }^{\alpha _{1}...\alpha
_{p}a_{1}...a_{r}}=\Gamma \left( \overset{\ast }{\tilde{\rho}},Id_{\overset{%
\ast }{E}}\right) \left( \overset{\ast }{\tilde{\delta}}_{\gamma }\right)
T_{\beta _{1}...\beta _{q}b_{1}...b_{s}}^{\alpha _{1}...\alpha
_{p}a_{1}...a_{r}} \\
\hspace*{8mm}+\left( \rho ,\eta \right) \overset{\ast }{H}_{\alpha \gamma
}^{\alpha _{1}}T_{\beta _{1}...\beta _{q}b_{1}...b_{s}}^{\alpha \alpha
_{2}...\alpha _{p}a_{1}...a_{r}}+...+\vspace*{2mm}\left( \rho ,\eta \right)
\overset{\ast }{H}_{\alpha \gamma }^{\alpha _{p}}T_{\beta _{1}...\beta
_{q}b_{1}...b_{s}}^{\alpha _{1}...\alpha _{p-1}\alpha a_{1}...a_{r}} \\
\hspace*{8mm}-\left( \rho ,\eta \right) \overset{\ast }{H}_{\beta _{1}\gamma
}^{\beta }T_{\beta \beta _{2}...\beta _{q}b_{1}...b_{s}}^{\alpha
_{1}...\alpha _{p}a_{1}...a_{r}}-...-\vspace*{2mm}\left( \rho ,\eta \right)
\overset{\ast }{H}_{\beta _{q}\gamma }^{\beta }T_{\beta _{1}...\beta
_{q-1}\beta b_{1}...b_{s}}^{\alpha _{1}...\alpha _{p}a_{1}...a_{r}} \\
\hspace*{8mm}-\left( \rho ,\eta \right) \overset{\ast }{H}_{a\gamma
}^{a_{1}}T_{\beta _{1}...\beta _{q}b_{1}...b_{s}}^{\alpha _{1}...\alpha
_{p}aa_{2}...a_{r}}-...-\vspace*{2mm}\left( \rho ,\eta \right) \overset{\ast
}{H}_{a\gamma }^{a_{r}}T_{\beta _{1}...\beta _{q}b_{1}...b_{s}}^{\alpha
_{1}...\alpha _{p}a_{1}...a_{r-1}a} \\
\hspace*{8mm}+\left( \rho ,\eta \right) \overset{\ast }{H}_{b_{1}\gamma
}^{b}T_{\beta _{1}...\beta _{q}bb_{2}...b_{s}}^{\alpha _{1}...\alpha
_{p}a_{1}...a_{r}}+\vspace*{2mm}...+\left( \rho ,\eta \right) \overset{\ast }%
{H}_{b_{s}\gamma }^{b}T_{\beta _{1}...\beta _{q}b_{1}...b_{s-1}b}^{\alpha
_{1}...\alpha _{p}a_{1}...a_{r}}%
\end{array}%
\leqno(4.6)
\end{equation*}%
\emph{and }%
\begin{equation*}
\begin{array}{l}
T_{\beta _{1}...\beta _{q}b_{1}...b_{s}}^{\alpha _{1}...\alpha
_{p}a_{1}...a_{r}}\mid ^{c}=\Gamma \left( \overset{\ast }{\tilde{\rho}},Id_{%
\overset{\ast }{E}}\right) \left( \overset{\cdot }{\tilde{\partial}}%
^{c}\right) T_{\beta _{1}...\beta _{q}b_{1}...b_{s}}^{\alpha _{1}...\alpha
_{p}a_{1}...a_{r}}\vspace*{2mm} \\
\hspace*{8mm}+\left( \rho ,\eta \right) \overset{\ast }{V}_{\alpha }^{\alpha
_{1}c}T_{\beta _{1}...\beta _{q}b_{1}...b_{s}}^{\alpha \alpha _{2}...\alpha
_{p}a_{1}...a_{r}}+...+\left( \rho ,\eta \right) \overset{\ast }{V}_{\alpha
}^{\alpha _{p}c}T_{\beta _{1}...\beta _{q}b_{1}...b_{s}}^{\alpha
_{1}...\alpha _{p-1}\alpha a_{1}...a_{r}}\vspace*{2mm} \\
\hspace*{8mm}-\left( \rho ,\eta \right) \overset{\ast }{V}_{\beta
_{1}}^{\beta c}T_{\beta \beta _{2}...\beta _{q}b_{1}...b_{s}}^{\alpha
_{1}...\alpha _{p}a_{1}...a_{r}}-...-\left( \rho ,\eta \right) \overset{\ast
}{V}_{\beta _{q}}^{\beta c}T_{\beta _{1}...\beta _{q-1}\beta
b_{1}...b_{s}}^{\alpha _{1}...\alpha _{p}a_{1}...a_{r}} \\
\hspace*{8mm}-\left( \rho ,\eta \right) \overset{\ast }{V}%
_{a}^{a_{1}c}T_{\beta _{1}...\beta _{q}b_{1}...b_{s}}^{\alpha _{1}...\alpha
_{p}aa_{2}...a_{r}}-...-\left( \rho ,\eta \right) \overset{\ast }{V}%
_{a}^{a_{r}c}T_{\beta _{1}...\beta _{q}b_{1}...b_{s}}^{\alpha _{1}...\alpha
_{p}a_{1}...a_{r-1}a}\vspace*{2mm} \\
\hspace*{8mm}+\left( \rho ,\eta \right) \overset{\ast }{V}%
_{b_{1}}^{bc}T_{\beta _{1}...\beta _{q}bb_{2}...b_{s}}^{\alpha _{1}...\alpha
_{p}a_{1}...a_{r}}...+\left( \rho ,\eta \right) \overset{\ast }{V}%
_{b_{s}}^{bc}T_{\beta _{1}...\beta _{q}b_{1}...b_{s-1}b}^{\alpha
_{1}...\alpha _{p}a_{1}...a_{r}}.%
\end{array}%
\leqno(4.7)
\end{equation*}

\textbf{Corollary 4.2 }\emph{In the particular case of Lie algebroids, }$%
\left( \eta ,h\right) =\left( Id_{M},Id_{M}\right) ,$\emph{\ we obtain }%
\begin{equation*}
\begin{array}{l}
T_{\beta _{1}...\beta _{q}b_{1}...b_{s}\mid \gamma }^{\alpha _{1}...\alpha
_{p}a_{1}...a_{r}}=\Gamma \left( \overset{\ast }{\tilde{\rho}},Id_{\overset{%
\ast }{E}}\right) \left( \overset{\ast }{\tilde{\delta}}_{\gamma }\right)
T_{\beta _{1}...\beta _{q}b_{1}...b_{s}}^{\alpha _{1}...\alpha
_{p}a_{1}...a_{r}} \\
\hspace*{8mm}+\rho \overset{\ast }{H}_{\alpha \gamma }^{\alpha _{1}}T_{\beta
_{1}...\beta _{q}b_{1}...b_{s}}^{\alpha \alpha _{2}...\alpha
_{p}a_{1}...a_{r}}+...+\vspace*{2mm}\rho \overset{\ast }{H}_{\alpha \gamma
}^{\alpha _{p}}T_{\beta _{1}...\beta _{q}b_{1}...b_{s}}^{\alpha
_{1}...\alpha _{p-1}\alpha a_{1}...a_{r}} \\
\hspace*{8mm}-\rho \overset{\ast }{H}_{\beta _{1}\gamma }^{\beta }T_{\beta
\beta _{2}...\beta _{q}b_{1}...b_{s}}^{\alpha _{1}...\alpha
_{p}a_{1}...a_{r}}-...-\vspace*{2mm}\rho \overset{\ast }{H}_{\beta
_{q}\gamma }^{\beta }T_{\beta _{1}...\beta _{q-1}\beta
b_{1}...b_{s}}^{\alpha _{1}...\alpha _{p}a_{1}...a_{r}} \\
\hspace*{8mm}-\rho \overset{\ast }{H}_{a\gamma }^{a_{1}}T_{\beta
_{1}...\beta _{q}b_{1}...b_{s}}^{\alpha _{1}...\alpha
_{p}aa_{2}...a_{r}}-...-\vspace*{2mm}\rho \overset{\ast }{H}_{a\gamma
}^{a_{r}}T_{\beta _{1}...\beta _{q}b_{1}...b_{s}}^{\alpha _{1}...\alpha
_{p}a_{1}...a_{r-1}a} \\
\hspace*{8mm}+\rho \overset{\ast }{H}_{b_{1}\gamma }^{b}T_{\beta
_{1}...\beta _{q}bb_{2}...b_{s}}^{\alpha _{1}...\alpha _{p}a_{1}...a_{r}}+%
\vspace*{2mm}...+\rho \overset{\ast }{H}_{b_{s}\gamma }^{b}T_{\beta
_{1}...\beta _{q}b_{1}...b_{s-1}b}^{\alpha _{1}...\alpha _{p}a_{1}...a_{r}}%
\end{array}%
\leqno(4.6)^{\prime }
\end{equation*}%
\emph{and }%
\begin{equation*}
\begin{array}{l}
T_{\beta _{1}...\beta _{q}b_{1}...b_{s}}^{\alpha _{1}...\alpha
_{p}a_{1}...a_{r}}\mid ^{c}=\Gamma \left( \overset{\ast }{\tilde{\rho}},Id_{%
\overset{\ast }{E}}\right) \left( \overset{\cdot }{\tilde{\partial}}%
^{c}\right) T_{\beta _{1}...\beta _{q}b_{1}...b_{s}}^{\alpha _{1}...\alpha
_{p}a_{1}...a_{r}} \\
\hspace*{8mm}+\rho \overset{\ast }{V}_{\alpha }^{\alpha _{1}c}T_{\beta
_{1}...\beta _{q}b_{1}...b_{s}}^{\alpha \alpha _{2}...\alpha
_{p}a_{1}...a_{r}}+...+\rho \overset{\ast }{V}_{\alpha }^{\alpha
_{p}c}T_{\beta _{1}...\beta _{q}b_{1}...b_{s}}^{\alpha _{1}...\alpha
_{p-1}\alpha a_{1}...a_{r}}\vspace*{2mm} \\
\hspace*{8mm}-\rho \overset{\ast }{V}_{\beta _{1}}^{\beta c}T_{\beta \beta
_{2}...\beta _{q}b_{1}...b_{s}}^{\alpha _{1}...\alpha
_{p}a_{1}...a_{r}}-...-\rho \overset{\ast }{V}_{\beta _{q}}^{\beta
c}T_{\beta _{1}...\beta _{q-1}\beta b_{1}...b_{s}}^{\alpha _{1}...\alpha
_{p}a_{1}...a_{r}} \\
\hspace*{8mm}-\rho \overset{\ast }{V}_{a}^{a_{1}c}T_{\beta _{1}...\beta
_{q}b_{1}...b_{s}}^{\alpha _{1}...\alpha _{p}aa_{2}...a_{r}}-...-\rho
\overset{\ast }{V}_{a}^{a_{r}c}T_{\beta _{1}...\beta
_{q}b_{1}...b_{s}}^{\alpha _{1}...\alpha _{p}a_{1}...a_{r-1}a}\vspace*{2mm}
\\
\hspace*{8mm}+\rho \overset{\ast }{V}_{b_{1}}^{bc}T_{\beta _{1}...\beta
_{q}bb_{2}...b_{s}}^{\alpha _{1}...\alpha _{p}a_{1}...a_{r}}...+\rho \overset%
{\ast }{V}_{b_{s}}^{bc}T_{\beta _{1}...\beta _{q}b_{1}...b_{s-1}b}^{\alpha
_{1}...\alpha _{p}a_{1}...a_{r}}.%
\end{array}%
\leqno(4.7)^{\prime }
\end{equation*}

\emph{In the classical case, }$\left( \rho ,\eta ,h\right) =\left(
Id_{TE},Id_{M},Id_{M}\right) ,$\emph{\ we obtain }%
\begin{equation*}
\begin{array}{l}
T_{j_{1}...j_{q}b_{1}...b_{s}\mid k}^{i_{1}...i_{p}a_{1}...a_{r}}=\vspace*{%
2mm}\overset{\ast }{\delta }_{k}\left(
T_{j_{1}...j_{q}b_{1}...b_{s}}^{i_{1}...i_{p}a_{1}...a_{r}}\right) \\
\hspace*{8mm}+\overset{\ast }{H}%
_{ik}^{i_{1}}T_{j_{1}...j_{q}b_{1}...b_{s}}^{ii_{2}...i_{p}a_{1}...a_{r}}+...+%
\vspace*{2mm}\overset{\ast }{H}_{ik}^{i_{p}}T_{\beta _{1}...\beta
_{q}b_{1}...b_{s}}^{i_{1}...i_{p-1}ia_{1}...a_{r}} \\
\hspace*{8mm}-\overset{\ast }{H}%
_{j_{1}k}^{j}T_{jj_{2}...j_{q}b_{1}...b_{s}}^{i_{1}...i_{p}a_{1}...a_{r}}-...-%
\vspace*{2mm}\overset{\ast }{H}%
_{j_{q}k}^{j}T_{j_{1}...j_{q-1}jb_{1}...b_{s}}^{\alpha _{1}...\alpha
_{p}a_{1}...a_{r}} \\
\hspace*{8mm}-\overset{\ast }{H}_{ak}^{a_{1}}T_{\beta _{1}...\beta
_{q}b_{1}...b_{s}}^{\alpha _{1}...\alpha _{p}aa_{2}...a_{r}}-...-\vspace*{2mm%
}\overset{\ast }{H}_{ak}^{a_{r}}T_{\beta _{1}...\beta
_{q}b_{1}...b_{s}}^{\alpha _{1}...\alpha _{p}a_{1}...a_{r-1}a} \\
\hspace*{8mm}+\overset{\ast }{H}_{b_{1}k}^{b}T_{\beta _{1}...\beta
_{q}bb_{2}...b_{s}}^{\alpha _{1}...\alpha _{p}a_{1}...a_{r}}+\vspace*{2mm}%
...+\overset{\ast }{H}_{b_{s}k}^{b}T_{\beta _{1}...\beta
_{q}b_{1}...b_{s-1}b}^{\alpha _{1}...\alpha _{p}a_{1}...a_{r}}%
\end{array}%
\leqno(4.6)^{\prime \prime }
\end{equation*}%
\emph{and }%
\begin{equation*}
\begin{array}{l}
T_{j_{1}...j_{q}b_{1}...b_{s}}^{i_{1}...i_{p}a_{1}...a_{r}}\mid ^{c}=\overset%
{\cdot }{\partial }^{c}\left(
T_{j_{1}...j_{q}b_{1}...b_{s}}^{i_{1}...i_{p}a_{1}...a_{r}}\right) \vspace*{%
2mm} \\
\hspace*{8mm}+\overset{\ast }{V}%
_{i}^{i_{1}c}T_{j_{1}...j_{q}b_{1}...b_{s}}^{ii_{2}...i_{p}a_{1}...a_{r}}+...+%
\overset{\ast }{V}%
_{i}^{i_{p}c}T_{j_{1}...j_{q}b_{1}...b_{s}}^{i_{1}...i_{p-1}ia_{1}...a_{r}}%
\vspace*{2mm} \\
\hspace*{8mm}-\overset{\ast }{V}%
_{j_{1}}^{jc}T_{jj_{2}...j_{q}b_{1}...b_{s}}^{i_{1}...i_{p}a_{1}...a_{r}}-...-%
\overset{\ast }{V}%
_{j_{q}}^{jc}T_{j_{1}...j_{q-1}jb_{1}...b_{s}}^{i_{1}...i_{p}a_{1}...a_{r}}
\\
\hspace*{8mm}-\overset{\ast }{V}%
_{a}^{a_{1}c}T_{j_{1}...j_{q}b_{1}...b_{s}}^{i_{1}...i_{p}aa_{2}...a_{r}}-...-%
\overset{\ast }{V}%
_{a}^{a_{r}c}T_{j_{1}...j_{q}b_{1}...b_{s}}^{i_{1}...i_{p}a_{1}...a_{r-1}a}%
\vspace*{2mm} \\
\hspace*{8mm}+\overset{\ast }{V}%
_{b_{1}}^{bc}T_{j_{1}...j_{q}bb_{2}...b_{s}}^{i_{1}...i_{p}a_{1}...a_{r}}...+%
\overset{\ast }{V}%
_{b_{s}}^{bc}T_{j_{1}...j_{q}b_{1}...b_{s-1}b}^{i_{1}...i_{p}a_{1}...a_{r}}.%
\end{array}%
\leqno(4.7)^{\prime \prime }
\end{equation*}

\textbf{Definition 4.2 }We assume that $\left( E,\pi ,M\right) =\left( F,\nu
,N\right) .$

If $\left( \rho ,\eta \right) \Gamma $ is a $\left( \rho ,\eta \right) $%
-connection for the vector bundle $\left( \overset{\ast }{E},\overset{\ast }{%
\pi },M\right) $ and
\begin{equation*}
\left( \left( \rho ,\eta \right) \overset{\ast }{H}_{bc}^{a},\left( \rho
,\eta \right) \overset{\ast }{\tilde{H}}_{bc}^{a},\left( \rho ,\eta \right)
\overset{\ast }{V}_{b}^{ac},\left( \rho ,\eta \right) \overset{\ast }{\tilde{%
V}}_{b}^{ac}\right)
\end{equation*}%
are the components of a distinguished linear $\left( \rho ,\eta \right) $%
\textit{-}connection for the generalized tangent bundle $\left( \left( \rho
,\eta \right) T\overset{\ast }{E},\left( \rho ,\eta \right) \tau _{\overset{%
\ast }{E}},\overset{\ast }{E}\right) $ such that
\begin{equation*}
\left( \rho ,\eta \right) \overset{\ast }{H}_{bc}^{a}=\left( \rho ,\eta
\right) \overset{\ast }{\tilde{H}}_{bc}^{a}\mbox{ and }\left( \rho ,\eta
\right) \overset{\ast }{V}_{b}^{ac}=\left( \rho ,\eta \right) \overset{\ast }%
{\tilde{V}}_{b}^{ac},
\end{equation*}%
then we will say that \emph{the generalized tangent bundle }$\!\left( \left(
\rho ,\eta \right) T\overset{\ast }{E},\left( \rho ,\eta \right) \tau _{%
\overset{\ast }{E}},\overset{\ast }{E}\right) $ \emph{is endowed with a
normal distinguished linear }$\left( \rho ,\eta \right) $\emph{-connection
on components }%
\begin{equation*}
\left( \left( \rho ,\eta \right) \overset{\ast }{H}_{bc}^{a},\left( \rho
,\eta \right) \overset{\ast }{V}_{b}^{ac}\right) .
\end{equation*}

In the particular case of Lie algebroids, $\left( \eta ,h\right) =\left(
Id_{M},Id_{M}\right) ,$\emph{\ }the components of a normal distinguished
linear $\left( \rho ,Id_{M}\right) $-connection $\left( \rho \overset{\ast }{%
H},\rho \overset{\ast }{V}\right) $ will be denoted $\left( \rho \overset{%
\ast }{H}_{bc}^{a},\rho \overset{\ast }{V}_{b}^{ac}\right) $.

In the classical case, $\left( \rho ,\eta ,h\right) =\left(
Id_{TE},Id_{M},Id_{M}\right) ,$\emph{\ }the components of a normal
distinguished linear $\left( Id_{TM},Id_{M}\right) $-connection $\left(
\overset{\ast }{H},\overset{\ast }{V}\right) $ will be denoted $\left(
\overset{\ast }{H}_{jk}^{i},\overset{\ast }{V}_{j}^{ik}\right) $.

\section{The $\left( \protect\rho ,\protect\eta \right) $%
-(pseudo)metrizability}

We consider the following diagram:
\begin{equation*}
\begin{array}{c}
\xymatrix{\overset{\ast }{E}\ar[d]_{\overset{\ast }{\pi }}&\left( F,\left[
,\right] _{F,h},\left( \rho ,\eta \right) \right)\ar[d]^\nu\\ M\ar[r]^h&N}%
\end{array}%
\end{equation*}%
where $\left( E,\pi ,M\right) \in \left\vert \mathbf{B}^{\mathbf{v}%
}\right\vert $ and $\left( \left( F,\nu ,M\right) ,\left[ ,\right]
_{F.h},\left( \rho ,\eta \right) \right) $ is a generalized Lie algebroid.
Let $\left( \rho ,\eta \right) \Gamma $ be a $\left( \rho ,\eta \right) $%
-connection for the vector bundle $\left( \overset{\ast }{E},\overset{\ast }{%
\pi },M\right) $ and let $\left( \left( \rho ,\eta \right) \overset{\ast }{H}%
,\left( \rho ,\eta \right) \overset{\ast }{V}\right) $ be a distinguished
linear $\left( \rho ,\eta \right) $-connection for the generalized tangent
bundle
\begin{equation*}
\left( \left( \rho ,\eta \right) T\overset{\ast }{E},\left( \rho ,\eta
\right) \tau _{\overset{\ast }{E}},\overset{\ast }{E}\right) .
\end{equation*}%
\noindent \textbf{Definition 5.1} A tensor $d$-field
\begin{equation*}
G=g_{\alpha \beta }d\tilde{z}^{\alpha }\otimes d\tilde{z}^{\beta
}+g^{ab}\delta \tilde{p}_{a}\otimes \delta \tilde{p}_{b}\in \mathcal{DT}%
_{20}^{02}\left( \left( \rho ,\eta \right) T\overset{\ast }{E},\left( \rho
,\eta \right) \tau _{\overset{\ast }{E}},\overset{\ast }{E}\right)
\end{equation*}%
will be called \emph{pseudometrical structure }if its components are
symmetric and the matrices $\left\Vert g_{\alpha \beta }\left( \overset{\ast
}{u}_{x}\right) \right\Vert $and $\left\Vert g_{ab}\left( \overset{\ast }{u}%
_{x}\right) \right\Vert $ are nondegenerate, for any point $\overset{\ast }{u%
}_{x}\in \overset{\ast }{E}.$

Moreover, if the matrices $\left\Vert g_{\alpha \beta }\left( \overset{\ast }%
{u}_{x}\right) \right\Vert $ and $\left\Vert g_{ab}\left( \overset{\ast }{u}%
_{x}\right) \right\Vert $ has constant signature, then the tensor $d$-field $%
G$ will be called \emph{metrical structure}\textit{.}

Let
\begin{equation*}
G=g_{\alpha \beta }d\tilde{z}^{\alpha }\otimes d\tilde{z}^{\beta
}+g^{ab}\delta \tilde{p}_{a}\otimes \delta \tilde{p}_{b}
\end{equation*}%
be a (pseudo)metrical structure. If $\alpha ,\beta \in \overline{1,p}$ and $%
a,b\in \overline{1,r},$ then for any vector local $\left( m+r\right) $-chart
$\left( U,\overset{\ast }{s}_{U}\right) $ of $\left( \overset{\ast }{E},%
\overset{\ast }{\pi },M\right) $, we consider the real functions
\begin{equation*}
\begin{array}{ccc}
\overset{\ast }{\pi }^{-1}\left( U\right) & ^{\underrightarrow{~\ \ \tilde{g}%
^{\beta \alpha }~\ \ }} & \mathbb{R}%
\end{array}%
\end{equation*}%
and
\begin{equation*}
\begin{array}{ccc}
\overset{\ast }{\pi }^{-1}\left( U\right) & ^{\underrightarrow{~\ \ \tilde{g}%
_{ba}~\ \ }} & \mathbb{R}%
\end{array}%
\end{equation*}%
such that%
\begin{equation*}
\begin{array}{c}
\left\Vert \tilde{g}^{\beta \alpha }\left( \overset{\ast }{u}_{x}\right)
\right\Vert =\left\Vert g_{\alpha \beta }\left( \overset{\ast }{u}%
_{x}\right) \right\Vert ^{-1}%
\end{array}%
\end{equation*}%
and
\begin{equation*}
\begin{array}{c}
\left\Vert \tilde{g}_{ba}\left( \overset{\ast }{u}_{x}\right) \right\Vert
=\left\Vert g^{ab}\left( \overset{\ast }{u}_{x}\right) \right\Vert ^{-1},%
\end{array}%
\end{equation*}%
for any $\overset{\ast }{u}_{x}\in \overset{\ast }{\pi }^{-1}\left( U\right)
\backslash \left\{ \overset{\ast }{0}_{x}\right\} $.

\bigskip \noindent \textbf{Definition 5.2} If around each point $x\in M$ it
exists a local vector $m+r$-chart $\left( U,\overset{\ast }{s}_{U}\right) $
and a local $m$-chart $\left( U,\xi _{U}\right) $ such that $g_{\alpha \beta
}\circ \overset{\ast }{s}_{U}^{-1}\circ \left( \xi _{U}^{-1}\times Id_{%
\mathbb{R}^{m}}\right) \left( x,p\right) $ and $g^{ab}\circ \overset{\ast }{s%
}_{U}^{-1}\circ \left( \xi _{U}^{-1}\times Id_{\mathbb{R}^{m}}\right) \left(
x,p\right) $ depends only on $x$, for any $\overset{\ast }{u}_{x}\in \overset%
{\ast }{\pi }^{-1}\left( U\right) ,$ then we will say that the \emph{%
(pseudo)metrical structure \ }%
\begin{equation*}
G=g_{\alpha \beta }d\tilde{z}^{\alpha }\otimes d\tilde{z}^{\beta
}+g^{ab}\delta \tilde{p}_{a}\otimes \delta \tilde{p}_{b}
\end{equation*}%
\emph{is a Riemannian (pseudo)metrical structure.}\textit{\ }

If only the condition is verified:

\textquotedblright $g_{\alpha \beta }\circ \overset{\ast }{s}_{U}^{-1}\circ
\left( \xi _{U}^{-1}\times Id_{\mathbb{R}^{m}}\right) \left( x,p\right) $%
\textit{\ depends only on }$x$\textit{, for any }$\overset{\ast }{u}_{x}\in
\overset{\ast }{\pi }^{-1}\left( U\right) $\textit{" res\-pec\-ti\-vely
\textquotedblright }$g^{ab}\circ \overset{\ast }{s}_{U}^{-1}\circ \left( \xi
_{U}^{-1}\times Id_{\mathbb{R}^{m}}\right) \left( x,p\right) $\textit{\
depends only on }$x$\textit{, for any }$\overset{\ast }{u}_{x}\in \overset{%
\ast }{\pi }^{-1}\left( U\right) $", then we will say that the \textit{%
(pseudo)metrical structure }$G$ \textit{is a }\emph{Riemannian }$\mathcal{H}$%
\emph{-(pseudo)metrical structure}\textit{\ }respectively a \emph{Riemannian
}$\mathcal{V}$\emph{-(pseudo)metrical structure.}

\noindent

\textbf{Definition 5.3} If around each point $x\in M$ there exists a local
vector $m+r$-chart $\left( U,\overset{\ast }{s}_{U}\right) $ and a local $m$%
-chart $\left( U,\xi _{U}\right) $ such that $g_{\alpha \beta }\circ \overset%
{\ast }{s}_{U}^{-1}\circ \left( \xi _{U}^{-1}\times Id_{\mathbb{R}%
^{m}}\right) \left( x,p\right) $ and $g^{ab}\circ \overset{\ast }{s}%
_{U}^{-1}\circ \left( \xi _{U}^{-1}\times Id_{\mathbb{R}^{m}}\right) \left(
x,p\right) $ depends only on~$p$\textit{, for any }$\overset{\ast }{u}%
_{x}\in \overset{\ast }{\pi }^{-1}\left( U\right) ,$ then we will say that
the \emph{(pseudo)metrical structure }%
\begin{equation*}
G=g_{\alpha \beta }d\tilde{z}^{\alpha }\otimes d\tilde{z}^{\beta
}+g^{ab}\delta \tilde{p}_{a}\otimes \delta \tilde{p}_{b}
\end{equation*}%
\emph{is a Minkowski (pseudo)metrical structure}$.$

If only the condition is verified:

\textquotedblright $g_{\alpha \beta }\circ \overset{\ast }{s}_{U}^{-1}\circ
\left( \xi _{U}^{-1}\times Id_{\mathbb{R}^{m}}\right) \left( x,p\right) $%
\textit{\ depends only on }$p$\textit{, for any }$\overset{\ast }{u}_{x}\in
\overset{\ast }{\pi }^{-1}\left( U\right) $" \textit{respectively
\textquotedblright }$g^{ab}\circ \overset{\ast }{s}_{U}^{-1}\circ \left( \xi
_{U}^{-1}\times Id_{\mathbb{R}^{m}}\right) \left( x,p\right) $ \textit{%
depends only on }$p$\textit{, for any }$\overset{\ast }{u}_{x}\in \overset{%
\ast }{\pi }^{-1}\left( U\right) $\textit{", }\noindent then we will say
that \emph{the (pseudo)metrical structure }$G$\emph{\ is a Minkowski }$%
\mathcal{H}$\emph{-(pseudo)metrical structure }respectively a \emph{%
Minkowski }$\mathcal{V}$\emph{-(pseudo)metrical structure.}

\textbf{Definition 5.4} If there exists a (pseudo)metrical structure%
\begin{equation*}
G=g_{\alpha \beta }d\tilde{z}^{\alpha }\otimes d\tilde{z}^{\beta
}+g^{ab}\delta \tilde{p}_{a}\otimes \delta \tilde{p}_{b}
\end{equation*}%
and a distinguished linear $\left( \rho ,\eta \right) $-connection%
\begin{equation*}
\left( \left( \rho ,\eta \right) \overset{\ast }{H},\left( \rho ,\eta
\right) \overset{\ast }{V}\right)
\end{equation*}%
such that\emph{\ }%
\begin{equation*}
\begin{array}{c}
\left( \rho ,\eta \right) \overset{\ast }{D}_{X}G=0,~\forall X\in \Gamma
\left( \left( \rho ,\eta \right) T\overset{\ast }{E},\left( \rho ,\eta
\right) \tau _{\overset{\ast }{E}},\overset{\ast }{E}\right) .%
\end{array}%
\leqno(5.1)
\end{equation*}%
then the generalized tangent bundle $\left( \left( \rho ,\eta \right) T%
\overset{\ast }{E},\left( \rho ,\eta \right) \tau _{\overset{\ast }{E}},%
\overset{\ast }{E}\right) $ will be called $(\rho ,\eta )$\emph{%
-(pseudo)metrizable}

Condition $\left( 5.1\right) $ is equivalent with the following equalities:
\begin{equation*}
\begin{array}{c}
g_{\alpha \beta \mid \gamma }=0,\,g_{~\ \ \mid \gamma }^{ab}=0,\,\,g_{\alpha
\beta }\mid ^{c}=0\,,\,\,g^{ab}\mid ^{c}=0.%
\end{array}%
\leqno(5.2)
\end{equation*}

If $g_{\alpha \beta \mid \gamma }{=}0$ and $\,g_{~\ \ \mid \gamma }^{ab}=0$,
then we will say that \emph{the vector bundle} $\left( \left( \rho ,\eta
\right) T\overset{\ast }{E},\left( \rho ,\eta \right) \tau _{\overset{\ast }{%
E}},\overset{\ast }{E}\right) $ \emph{is }$\mathcal{H}$\emph{-}$(\rho ,\eta
) $\emph{-(pseudo)metrizable.}

If $g_{\alpha \beta }|^{c}{=}0$ and $\,g^{ab}\mid ^{c}=0$, then we will say
that \emph{the vector bundle} $\left( \left( \rho ,\eta \right) T\overset{%
\ast }{E},\left( \rho ,\eta \right) \tau _{\overset{\ast }{E}},\overset{\ast
}{E}\right) $ \emph{is }$\mathcal{V}$\emph{-}$(\rho ,\eta )$\emph{%
-(pseudo)\-metrizable.}

\bigskip \noindent \textbf{Theorem 5.1} \emph{If} $\left( \left( \rho ,\eta
\right) \overset{\ast }{\mathring{H}},\left( \rho ,\eta \right) \overset{%
\ast }{\mathring{V}}\right) $ \emph{is a distinguished linear }$\left( \rho
,\eta \right) $\emph{-connection for the generalized tangent bundle }$\left(
\left( \rho ,\eta \right) T\overset{\ast }{E},\left( \rho ,\eta \right) \tau
_{\overset{\ast }{E}},\overset{\ast }{E}\right) $ \emph{\ and }$G=g_{\alpha
\beta }d\tilde{z}^{\alpha }\otimes d\tilde{z}^{\beta }+g^{ab}\delta \tilde{p}%
_{a}\otimes \delta \tilde{p}_{b}$ \emph{is a (pseudo)metrical structure,
then the following real local functions:}%
\begin{equation*}
(5.3)%
\begin{array}{ll}
\left( \rho ,\eta \right) \overset{\ast }{H}_{\beta \gamma }^{\alpha }\!\! &
=\displaystyle\frac{1}{2}\tilde{g}^{\alpha \varepsilon }\left( \Gamma \left(
\overset{\ast }{\tilde{\rho}},Id_{E}\right) \left( \overset{\ast }{\tilde{%
\delta}}_{\gamma }\right) g_{\varepsilon \beta }+\Gamma \left( \overset{\ast
}{\tilde{\rho}},Id_{E}\right) \left( \overset{\ast }{\tilde{\delta}}_{\beta
}\right) g_{\varepsilon \gamma }-\Gamma \left( \overset{\ast }{\tilde{\rho}}%
,Id_{E}\right) \left( \overset{\ast }{\tilde{\delta}}_{\varepsilon }\right)
g_{\beta \gamma }\right. \vspace*{1mm} \\
& \left. +g_{\theta \varepsilon }L_{\gamma \beta }^{\theta }\circ h\circ
\overset{\ast }{\pi }-g_{\beta \theta }L_{\gamma \varepsilon }^{\theta
}\circ h\circ \overset{\ast }{\pi }-g_{\theta \gamma }L_{\beta \varepsilon
}^{\theta }\circ h\circ \overset{\ast }{\pi }\right) ,\vspace*{2mm} \\
\left( \rho ,\eta \right) \overset{\ast }{H}_{b\gamma }^{a}\!\! & =\left(
\rho ,\eta \right) \overset{\ast }{\mathring{H}}_{b\gamma }^{a}+\displaystyle%
\frac{1}{2}\tilde{g}_{be}g_{~\ \ \ \overset{0}{\mid }\gamma }^{ea},\vspace*{%
2mm} \\
\left( \rho ,\eta \right) \overset{\ast }{V}_{\beta }^{\alpha c}\!\! &
=\left( \rho ,\eta \right) \overset{\ast }{\mathring{V}}_{\beta }^{\alpha c}+%
\displaystyle\frac{1}{2}\tilde{g}^{\alpha \varepsilon }g_{\varepsilon \beta }%
\overset{0}{\mid }^{c},\vspace*{2mm} \\
\left( \rho ,\eta \right) \overset{\ast }{V}_{b}^{ac}\!\! & =\displaystyle%
\frac{1}{2}\tilde{g}_{be}\left( \Gamma \left( \overset{\ast }{\tilde{\rho}}%
,Id_{E}\right) \left( \overset{\cdot }{\tilde{\partial}}^{c}\right)
g^{ea}+\Gamma \left( \overset{\ast }{\tilde{\rho}},Id_{E}\right) \left(
\overset{\cdot }{\tilde{\partial}}^{a}\right) g^{ec}-\Gamma \left( \overset{%
\ast }{\tilde{\rho}},Id_{E}\right) \left( \overset{\cdot }{\tilde{\partial}}%
^{e}\right) g^{ac}\right)%
\end{array}%
\end{equation*}%
\emph{are components of a distinguished linear }$\left( \rho ,\eta \right) $%
\emph{-connection such that the generalized tangent bundle }$\left( \left(
\rho ,\eta \right) T\overset{\ast }{E},\left( \rho ,\eta \right) \tau _{%
\overset{\ast }{E}},\overset{\ast }{E}\right) $ \emph{becomes }$\left( \rho
,\eta \right) $\emph{-(pseudo)metrizable.}

\textbf{Corollary 5.1 }\emph{In the particular case of Lie algebroids, }$%
\left( \eta ,h\right) =\left( Id_{M},Id_{M}\right) ,$\emph{\ then we obtain}
\begin{equation*}
(5.3)^{\prime }%
\begin{array}{ll}
\rho \overset{\ast }{H}_{\beta \gamma }^{\alpha }\!\! & =\displaystyle\frac{1%
}{2}\tilde{g}^{\alpha \varepsilon }\left( \Gamma \left( \overset{\ast }{%
\tilde{\rho}},Id_{E}\right) \left( \overset{\ast }{\tilde{\delta}}_{\gamma
}\right) g_{\varepsilon \beta }+\Gamma \left( \overset{\ast }{\tilde{\rho}}%
,Id_{E}\right) \left( \overset{\ast }{\tilde{\delta}}_{\beta }\right)
g_{\varepsilon \gamma }-\Gamma \left( \overset{\ast }{\tilde{\rho}}%
,Id_{E}\right) \left( \overset{\ast }{\tilde{\delta}}_{\varepsilon }\right)
g_{\beta \gamma }\right. \vspace*{1mm} \\
& \left. +g_{\theta \varepsilon }L_{\gamma \beta }^{\theta }\circ h\circ
\overset{\ast }{\pi }-g_{\beta \theta }L_{\gamma \varepsilon }^{\theta
}\circ h\circ \overset{\ast }{\pi }-g_{\theta \gamma }L_{\beta \varepsilon
}^{\theta }\circ h\circ \overset{\ast }{\pi }\right) ,\vspace*{2mm} \\
\rho \overset{\ast }{H}_{b\gamma }^{a}\!\! & =\rho \overset{\ast }{\mathring{%
H}}_{b\gamma }^{a}+\displaystyle\frac{1}{2}\tilde{g}_{be}g_{~\ \ \ \overset{0%
}{\mid }\gamma }^{ea},\vspace*{2mm} \\
\rho \overset{\ast }{V}_{\beta }^{\alpha c}\!\! & =\rho \overset{\ast }{%
\mathring{V}}_{\beta }^{\alpha c}+\displaystyle\frac{1}{2}\tilde{g}^{\alpha
\varepsilon }g_{\varepsilon \beta }\overset{0}{\mid }^{c},\vspace*{2mm} \\
\rho \overset{\ast }{V}_{b}^{ac}\!\! & =\displaystyle\frac{1}{2}\tilde{g}%
_{be}\left( \Gamma \left( \overset{\ast }{\tilde{\rho}},Id_{E}\right) \left(
\overset{\cdot }{\tilde{\partial}}^{c}\right) g^{ea}+\Gamma \left( \overset{%
\ast }{\tilde{\rho}},Id_{E}\right) \left( \overset{\cdot }{\tilde{\partial}}%
^{a}\right) g^{ec}-\Gamma \left( \overset{\ast }{\tilde{\rho}},Id_{E}\right)
\left( \overset{\cdot }{\tilde{\partial}}^{e}\right) g^{ac}\right)%
\end{array}%
\end{equation*}

\emph{In the classicale case, }$\left( \rho ,\eta ,h\right) =\left(
Id_{TE},Id_{M},Id_{M}\right) ,$\emph{\ then we obtain}
\begin{equation*}
\begin{array}{ll}
\overset{\ast }{H}_{jk}^{i}\!\! & =\displaystyle\frac{1}{2}\tilde{g}%
^{ih}\left( \overset{\ast }{\delta }_{k}g_{hj}+\overset{\ast }{\delta }%
_{j}g_{hk}-\overset{\ast }{\delta }_{h}g_{jk}\right) \\
\overset{\ast }{H}_{bk}^{a}\!\! & =\overset{\ast }{\mathring{H}}_{bk}^{a}+%
\displaystyle\frac{1}{2}\tilde{g}_{be}g_{~\ \ \ \overset{0}{\mid }k}^{ea},%
\vspace*{2mm} \\
\overset{\ast }{V}_{j}^{ic}\!\! & =\overset{\ast }{\mathring{V}}_{j}^{ic}+%
\displaystyle\frac{1}{2}\tilde{g}^{ih}g_{hj}\overset{0}{\mid }^{c},\vspace*{%
2mm} \\
\overset{\ast }{V}_{b}^{ac}\!\! & =\displaystyle\frac{1}{2}\tilde{g}%
_{be}\left( \overset{\cdot }{\partial }^{c}g^{ea}+\overset{\cdot }{\partial }%
^{a}g^{ec}-\overset{\cdot }{\partial }^{e}g^{ac}\right)%
\end{array}%
\leqno(5.3)^{\prime \prime }
\end{equation*}

\bigskip \noindent \textbf{Theorem 5.2 }\emph{If the distinguished linear }$%
\left( \rho ,\eta \right) $\emph{-connection} $\left( \left( \rho ,\eta
\right) \overset{0}{H},\left( \rho ,\eta \right) \overset{0}{V}\right) $
\emph{coincides with the Berwald linear }$\left( \rho ,\eta \right) $\emph{%
-connection in the previous theorem, then the local real functions: }
\begin{equation*}
(5.4)%
\begin{array}{ll}
\left( \rho ,\eta \right) \overset{c}{\overset{\ast }{H}}_{\beta \gamma
}^{\alpha }\!\!\! & =\displaystyle\frac{1}{2}\tilde{g}^{\alpha \varepsilon
}\left( \Gamma \left( \overset{\ast }{\tilde{\rho}},Id_{E}\right) \left(
\overset{\ast }{\tilde{\delta}}_{\gamma }\right) g_{\varepsilon \beta
}+\Gamma \left( \overset{\ast }{\tilde{\rho}},Id_{E}\right) \left( \overset{%
\ast }{\tilde{\delta}}_{\beta }\right) g_{\varepsilon \gamma }\right.
\vspace*{1mm} \\
& -\Gamma \left( \overset{\ast }{\tilde{\rho}},Id_{E}\right) \left( \overset{%
\ast }{\tilde{\delta}}_{\varepsilon }\right) g_{\beta \gamma }+g_{\theta
\varepsilon }L_{\gamma \beta }^{\theta }\circ h\circ \overset{\ast }{\pi }%
\vspace*{1mm}\left. -g_{\beta \theta }L_{\gamma \varepsilon }^{\theta }\circ
h\circ \overset{\ast }{\pi }\vspace*{1mm}-g_{\theta \gamma }L_{\beta
\varepsilon }^{\theta }\circ h\circ \overset{\ast }{\pi }\vspace*{1mm}%
\right) , \\
\left( \rho ,\eta \right) \overset{c}{\overset{\ast }{H}}_{b\gamma
}^{a}\!\!\! & =\displaystyle\frac{\partial \left( \rho ,\eta \right) \Gamma
_{b\gamma }}{\partial p_{a}}+\frac{1}{2}\tilde{g}_{be}g_{~\ \ \ \overset{0}{%
\mid }\gamma }^{ea},\vspace*{2mm} \\
\left( \rho ,\eta \right) \overset{c}{\overset{\ast }{V}}_{\beta }^{\alpha
c}\!\! & =\displaystyle\frac{1}{2}\tilde{g}_{\beta \varepsilon }\frac{%
\partial g^{\varepsilon \alpha }}{\partial p_{c}},\vspace*{2mm} \\
\left( \rho ,\eta \right) \overset{c}{\overset{\ast }{V}}_{b}^{ac}\!\!\! & =%
\displaystyle\frac{1}{2}\tilde{g}_{be}\left( \frac{\partial g^{ea}}{\partial
p_{c}}+\frac{\partial g^{ec}}{\partial p_{a}}-\frac{\partial g^{ac}}{%
\partial p_{e}}\right)%
\end{array}%
\hspace*{-6mm}
\end{equation*}%
\emph{are the components of a distinguished linear }$\left( \rho ,\eta
\right) $\emph{-connection such that the generalized tangent bundle }$\left(
\left( \rho ,\eta \right) T\overset{\ast }{E},\left( \rho ,\eta \right) \tau
_{\overset{\ast }{E}},\overset{\ast }{E}\right) $ \emph{becomes }$\left(
\rho ,\eta \right) $\emph{-(pseudo)metrizable.}

\emph{Moreover, if the (pseudo)metrical structure }$G$\emph{\ is }$\mathcal{H%
}$\emph{- and }$\mathcal{V}$\emph{-Rieman\-nian, then the local real
functions:}\smallskip \noindent
\begin{equation*}
\begin{array}{cl}
(\rho ,\eta )\overset{c}{\overset{\ast }{H}}_{\beta \gamma }^{\alpha } &
\overset{c}{\overset{\ast }{H}}{=}\frac{1}{2}\tilde{g}^{\alpha \varepsilon
}\left( \rho _{\gamma }^{k}{\circ }h{\circ }\overset{\ast }{\pi }\frac{%
\partial g_{\varepsilon \beta }}{\partial x^{k}}+\rho _{\beta }^{j}{\circ }h{%
\circ }\overset{\ast }{\pi }\frac{\partial g_{\varepsilon \gamma }}{\partial
x^{j}}-\rho _{\varepsilon }^{e}{\circ }h{\circ }\overset{\ast }{\pi }\frac{%
\partial g_{\beta \gamma }}{\partial x^{e}}+\right. \\
& \left. +g_{\theta \varepsilon }L_{\gamma \beta }^{\theta }{\circ }h{\circ }%
\overset{\ast }{\pi }-g_{\beta \theta }L_{\gamma \varepsilon }^{\theta }{%
\circ }h{\circ }\overset{\ast }{\pi }-g_{\theta \gamma }L_{\beta \varepsilon
}^{\theta }{\circ }h{\circ }\overset{\ast }{\pi }\right) ,\vspace*{1mm} \\
\left( \rho ,\eta \right) \overset{c}{\overset{\ast }{H}}_{b\gamma }^{a} & {=%
}\frac{\partial \left( \rho ,\eta \right) \Gamma _{b\gamma }}{\partial p_{a}}%
+\frac{1}{2}\tilde{g}_{be}\left( \rho _{\gamma }^{i}{\circ }h{\circ }\overset%
{\ast }{\pi }\frac{\partial g^{ea}}{\partial x^{i}}-\frac{\partial \left(
\rho ,\eta \right) \Gamma _{d\gamma }}{\partial p_{e}}g^{da}-\frac{\partial
\left( \rho ,\eta \right) \Gamma _{d\gamma }}{\partial p_{a}}g^{ed}\right) ,
\\
\left( \rho ,\eta \right) \overset{c}{\overset{\ast }{V}}_{\beta }^{\alpha c}
& =0, \\
\left( \rho ,\eta \right) \overset{c}{\overset{\ast }{V}}_{b}^{ac} & =0.%
\end{array}%
\leqno(5.5)
\end{equation*}%
\emph{are the components of a distinguished linear }$\left( \rho ,\eta
\right) $\emph{-connection such that the generalized tangent bundle }$\left(
\left( \rho ,\eta \right) T\overset{\ast }{E},\left( \rho ,\eta \right) \tau
_{\overset{\ast }{E}},\overset{\ast }{E}\right) $ \emph{becomes }$\left(
\rho ,\eta \right) $\emph{-(pseudo)metrizable.}

\textbf{Corollary 5.2 }\emph{In the particular case of Lie algebroids, }$%
\left( \eta ,h\right) =\left( Id_{M},Id_{M}\right) ,$\emph{\ then we obtain}
\emph{\ }
\begin{equation*}
\begin{array}{ll}
\rho \overset{c}{\overset{\ast }{H}}_{\beta \gamma }^{\alpha }\!\!\! & =%
\displaystyle\frac{1}{2}\tilde{g}^{\alpha \varepsilon }\left( \Gamma \left(
\overset{\ast }{\tilde{\rho}},Id_{E}\right) \left( \overset{\ast }{\tilde{%
\delta}}_{\gamma }\right) g_{\varepsilon \beta }+\Gamma \left( \overset{\ast
}{\tilde{\rho}},Id_{E}\right) \left( \overset{\ast }{\tilde{\delta}}_{\beta
}\right) g_{\varepsilon \gamma }\right. \vspace*{1mm} \\
& -\Gamma \left( \overset{\ast }{\tilde{\rho}},Id_{E}\right) \left( \overset{%
\ast }{\tilde{\delta}}_{\varepsilon }\right) g_{\beta \gamma }+g_{\theta
\varepsilon }L_{\gamma \beta }^{\theta }\circ \overset{\ast }{\pi }\vspace*{%
1mm}\left. -g_{\beta \theta }L_{\gamma \varepsilon }^{\theta }\circ \overset{%
\ast }{\pi }\vspace*{1mm}-g_{\theta \gamma }L_{\beta \varepsilon }^{\theta
}\circ \overset{\ast }{\pi }\vspace*{1mm}\right) , \\
\rho \overset{c}{\overset{\ast }{H}}_{b\gamma }^{a}\!\!\! & =\displaystyle%
\frac{\partial \rho \Gamma _{b\gamma }}{\partial p_{a}}+\frac{1}{2}\tilde{g}%
_{be}g_{~\ \ \ \overset{0}{\mid }\gamma }^{ea},\vspace*{2mm} \\
\rho \overset{c}{\overset{\ast }{V}}_{\beta }^{\alpha c}\!\! & =\displaystyle%
\frac{1}{2}\tilde{g}_{\beta \varepsilon }\frac{\partial g^{\varepsilon
\alpha }}{\partial p_{c}},\vspace*{2mm} \\
\rho \overset{c}{\overset{\ast }{V}}_{b}^{ac}\!\!\! & =\displaystyle\frac{1}{%
2}\tilde{g}_{be}\left( \frac{\partial g^{ea}}{\partial p_{c}}+\frac{\partial
g^{ec}}{\partial p_{a}}-\frac{\partial g^{ac}}{\partial p_{e}}\right)%
\end{array}%
\hspace*{-6mm}\leqno(5.4)^{\prime }
\end{equation*}

\emph{If the (pseudo)metrical structure }$G$\emph{\ is }$\mathcal{H}$\emph{-
and }$\mathcal{V}$\emph{-Rieman\-nian, then }%
\begin{equation*}
\begin{array}{l}
\begin{array}{ll}
\rho \overset{c}{\overset{\ast }{H}}_{\beta \gamma }^{\alpha }\!\!\! & =%
\displaystyle\frac{1}{2}\tilde{g}^{\alpha \varepsilon }\left( \Gamma \left(
\overset{\ast }{\tilde{\rho}},Id_{E}\right) \left( \overset{\ast }{\tilde{%
\delta}}_{\gamma }\right) g_{\varepsilon \beta }+\Gamma \left( \overset{\ast
}{\tilde{\rho}},Id_{E}\right) \left( \overset{\ast }{\tilde{\delta}}_{\beta
}\right) g_{\varepsilon \gamma }\right. \vspace*{1mm} \\
& -\Gamma \left( \overset{\ast }{\tilde{\rho}},Id_{E}\right) \left( \overset{%
\ast }{\tilde{\delta}}_{\varepsilon }\right) g_{\beta \gamma }+g_{\theta
\varepsilon }L_{\gamma \beta }^{\theta }\circ \overset{\ast }{\pi }\vspace*{%
1mm}\left. -g_{\beta \theta }L_{\gamma \varepsilon }^{\theta }\circ \overset{%
\ast }{\pi }\vspace*{1mm}-g_{\theta \gamma }L_{\beta \varepsilon }^{\theta
}\circ \overset{\ast }{\pi }\vspace*{1mm}\right) , \\
\rho \overset{c}{\overset{\ast }{H}}_{b\gamma }^{a}\!\!\! & =\displaystyle%
\frac{\partial \rho \Gamma _{b\gamma }}{\partial p_{a}}+\frac{1}{2}\tilde{g}%
_{be}g_{~\ \ \ \overset{0}{\mid }\gamma }^{ea},\vspace*{2mm} \\
\rho \overset{c}{\overset{\ast }{V}}_{\beta }^{\alpha c} & =0 \\
\rho \overset{c}{\overset{\ast }{V}}_{b}^{ac} & =0%
\end{array}%
\end{array}%
\leqno(5.5)^{\prime }
\end{equation*}

\emph{In the classicale case, }$\left( \rho ,\eta ,h\right) =\left(
Id_{TE},Id_{M},Id_{M}\right) ,$\emph{\ then we obtain}%
\begin{equation*}
\begin{array}{ll}
\overset{c}{\overset{\ast }{H}}_{jk}^{i}\!\!\! & =\displaystyle\frac{1}{2}%
\tilde{g}^{ih}\left( \overset{\ast }{\delta }_{k}g_{hj}+\overset{\ast }{%
\delta }_{j}g_{hk}-\overset{\ast }{\tilde{\delta}}_{h}g_{jk}\right) \\
\overset{c}{\overset{\ast }{H}}_{bk}^{a}\!\!\! & =\displaystyle\frac{%
\partial \Gamma _{bk}}{\partial p_{a}}+\frac{1}{2}\tilde{g}_{be}g_{~\ \ \
\overset{0}{\mid }k}^{ea},\vspace*{2mm} \\
\overset{c}{\overset{\ast }{V}}_{j}^{ic}\!\! & =\displaystyle\frac{1}{2}%
\tilde{g}_{jh}\frac{\partial g^{hi}}{\partial p_{c}},\vspace*{2mm} \\
\overset{c}{\overset{\ast }{V}}_{b}^{ac}\!\!\! & =\displaystyle\frac{1}{2}%
\tilde{g}_{be}\left( \frac{\partial g^{ea}}{\partial p_{c}}+\frac{\partial
g^{ec}}{\partial p_{a}}-\frac{\partial g^{ac}}{\partial p_{e}}\right)%
\end{array}%
\leqno(5.4)^{\prime \prime }
\end{equation*}

\emph{If the (pseudo)metrical structure }$G$\emph{\ is }$\mathcal{H}$\emph{-
and }$\mathcal{V}$\emph{-Rieman\-nian, then }%
\begin{equation*}
\begin{array}{l}
\overset{c}{\overset{\ast }{H}}_{jk}^{i}{=}\displaystyle\frac{1}{2}\tilde{g}%
^{ih}\left( \frac{\partial g_{hj}}{\partial x^{k}}+\frac{\partial g_{hk}}{%
\partial x^{j}}-\frac{\partial g_{jk}}{\partial x^{h}}\right)  \\
\overset{c}{\overset{\ast }{H}}_{bk}^{a}{=}\frac{\partial \Gamma _{bk}}{%
\partial p_{a}}+\frac{1}{2}\tilde{g}_{be}\left( \frac{\partial g^{ea}}{%
\partial x^{i}}-\frac{\partial \Gamma _{dk}}{\partial p_{e}}g^{da}-\frac{%
\partial \Gamma _{dk}}{\partial p_{a}}g^{ed}\right) ,\vspace*{2mm} \\
\overset{c}{\overset{\ast }{V}}_{j}^{ic}=0,\ \overset{c}{\overset{\ast }{V}}%
_{b}^{ac}=0.%
\end{array}%
\leqno(5.5)^{\prime \prime }
\end{equation*}%
\textbf{\ }

\bigskip \noindent \textbf{Theorem 5.3} \emph{Let }$\left( \rho ,\eta
\right) \Gamma $\emph{\ be a }$\left( \rho ,\eta \right) $\emph{-connection
for the vector bundle }$\left( \overset{\ast }{E},\overset{\ast }{\pi }%
,M\right) .$ \emph{Let }%
\begin{equation*}
\left( \left( \rho ,\eta \right) \overset{\ast }{\mathring{H}},\left( \rho
,\eta \right) \overset{\ast }{\mathring{V}}\right)
\end{equation*}%
\emph{be a distinguished linear }$\left( \rho ,\eta \right) $\emph{%
-connection for }%
\begin{equation*}
\left( \left( \rho ,\eta \right) T\overset{\ast }{E},\left( \rho ,\eta
\right) \tau _{\overset{\ast }{E}},\overset{\ast }{E}\right)
\end{equation*}
\emph{and let }%
\begin{equation*}
G=g_{\alpha \beta }d\tilde{z}^{\alpha }\otimes d\tilde{z}^{\beta
}+g^{ab}\delta \tilde{p}_{a}\otimes \delta \tilde{p}_{b}
\end{equation*}%
\emph{be a (pseudo)metrical structure. Let }
\begin{equation*}
\begin{array}{ll}
O_{\beta \gamma }^{\alpha \varepsilon }=\frac{1}{2}\left( \delta _{\beta
}^{\alpha }\delta _{\gamma }^{\varepsilon }-g_{\beta \gamma }\tilde{g}%
^{\alpha \varepsilon }\right) , & O_{\beta \gamma }^{\ast \alpha \varepsilon
}=\displaystyle\frac{1}{2}\left( \delta _{\beta }^{\alpha }\delta _{\gamma
}^{\varepsilon }+g_{\beta \gamma }\tilde{g}^{\alpha \varepsilon }\right) ,%
\vspace*{2mm} \\
O_{bc}^{ae}=\displaystyle\frac{1}{2}\left( \delta _{b}^{a}\delta _{c}^{e}-%
\tilde{g}_{bc}g^{ae}\right) , & O_{bc}^{\ast ae}=\frac{1}{2}\left( \delta
_{b}^{a}\delta _{c}^{e}+\tilde{g}_{bc}g^{ae}\right) ,%
\end{array}%
\leqno(5.6)
\end{equation*}%
\emph{be the Obata operators}.

\emph{If the real local functions }$X_{\beta \gamma }^{\alpha },X_{\beta
}^{\alpha c},Y_{b\gamma }^{a},Y_{b}^{ac}$ \emph{are components of tensor
fields,} \emph{then the local real functions given in the following: }%
\vspace*{-3mm}
\begin{equation*}
\begin{array}{ll}
\left( \rho ,\eta \right) \overset{\ast }{H}_{\beta \gamma }^{\alpha }%
\vspace*{1mm} & =\left( \rho ,\eta \right) \overset{c}{\overset{\ast }{H}}%
_{\beta \gamma }^{\alpha }+O_{\gamma \eta }^{\alpha \varepsilon
}X_{\varepsilon \beta }^{\eta }, \\
\left( \rho ,\eta \right) \overset{\ast }{H}_{b\gamma }^{a}\vspace*{1mm} &
=\left( \rho ,\eta \right) \overset{c}{\overset{\ast }{H}}_{b\gamma
}^{a}+O_{bd}^{ae}Y_{e\gamma }^{d},\vspace*{1mm} \\
\left( \rho ,\eta \right) \overset{\ast }{V}_{\beta }^{\alpha c} & =\left(
\rho ,\eta \right) \overset{c}{\overset{\ast }{V}}_{\beta }^{\alpha
c}+O_{\beta \eta }^{\ast \alpha \varepsilon }X_{\varepsilon }^{\eta c},%
\vspace*{1mm} \\
\left( \rho ,\eta \right) \overset{\ast }{V}_{b}^{ac} & =\left( \rho ,\eta
\right) \overset{c}{\overset{\ast }{V}}_{b}^{ac}+O_{bd}^{\ast ae}Y_{e}^{dc},%
\end{array}%
\leqno(5.7)
\end{equation*}%
\emph{are the components of a distinguished linear }$\left( \rho ,\eta
\right) $\emph{-connection such that the generalized tangent bundle }$\left(
\left( \rho ,\eta \right) T\overset{\ast }{E},\left( \rho ,\eta \right) \tau
_{\overset{\ast }{E}},\overset{\ast }{E}\right) $ \emph{becomes }$\left(
\rho ,\eta \right) $\emph{-(pseudo)metrizable.}

\textbf{Corollary 5.2 }\emph{In the particular case of Lie algebroids, }$%
\left( \eta ,h\right) =\left( Id_{M},Id_{M}\right) ,$\emph{\ then we obtain}%
\begin{equation*}
\begin{array}{ll}
\rho \overset{\ast }{H}_{\beta \gamma }^{\alpha }\vspace*{1mm} & =\rho
\overset{c}{\overset{\ast }{H}}_{\beta \gamma }^{\alpha }+O_{\gamma \eta
}^{\alpha \varepsilon }X_{\varepsilon \beta }^{\eta }, \\
\rho \overset{\ast }{H}_{b\gamma }^{a}\vspace*{1mm} & =\rho \overset{c}{%
\overset{\ast }{H}}_{b\gamma }^{a}+O_{bd}^{ae}Y_{e\gamma }^{d},\vspace*{1mm}
\\
\rho \overset{\ast }{V}_{\beta }^{\alpha c} & =\rho \overset{c}{\overset{%
\ast }{V}}_{\beta }^{\alpha c}+O_{\beta \eta }^{\ast \alpha \varepsilon
}X_{\varepsilon }^{\eta c},\vspace*{1mm} \\
\rho \overset{\ast }{V}_{b}^{ac} & =\rho \overset{c}{\overset{\ast }{V}}%
_{b}^{ac}+O_{bd}^{\ast ae}Y_{e}^{dc},%
\end{array}%
\leqno(5.7)^{\prime }
\end{equation*}

\emph{In the classicale case, }$\left( \rho ,\eta ,h\right) =\left(
Id_{TE},Id_{M},Id_{M}\right) ,$\emph{\ then we obtain }%
\begin{equation*}
\begin{array}{ll}
\overset{\ast }{H}_{jk}^{i}\vspace*{1mm} & =\overset{c}{\overset{\ast }{H}}%
_{jk}^{i}+O_{kl}^{ih}X_{hj}^{l}, \\
\overset{\ast }{H}_{bk}^{a}\vspace*{1mm} & =\overset{c}{\overset{\ast }{H}}%
_{bk}^{a}+O_{bd}^{ae}Y_{ek}^{d},\vspace*{1mm} \\
\overset{\ast }{V}_{j}^{ic} & =\overset{c}{\overset{\ast }{V}}%
_{j}^{ic}+O_{jl}^{\ast ih}X_{h}^{lc},\vspace*{1mm} \\
\overset{\ast }{V}_{b}^{ac} & =\overset{c}{\overset{\ast }{V}}%
_{b}^{ac}+O_{bd}^{\ast ae}Y_{e}^{dc},%
\end{array}%
\leqno(5.7)^{\prime \prime }
\end{equation*}

\bigskip \noindent \textbf{Theorem 5.3 }\emph{Let }$\left( \rho ,\eta
\right) \Gamma $\emph{\ be a }$\left( \rho ,\eta \right) $\emph{-connection
for the vector bundle }$\left( \overset{\ast }{E},\overset{\ast }{\pi }%
,M\right) .$ \emph{If }%
\begin{equation*}
\left( \left( \rho ,\eta \right) \overset{\ast }{\mathring{H}},\left( \rho
,\eta \right) \overset{\ast }{\mathring{V}}\right)
\end{equation*}%
\emph{is a distinguished linear }$\left( \rho ,\eta \right) $\emph{%
-connection for the generalized tangent bundle }$\left( \left( \rho ,\eta
\right) T\overset{\ast }{E},\left( \rho ,\eta \right) \tau _{\overset{\ast }{%
E}},\overset{\ast }{E}\right) $\emph{\ and }%
\begin{equation*}
G=g_{\alpha \beta }d\tilde{z}^{\alpha }\otimes d\tilde{z}^{\beta
}+g^{ab}\delta \tilde{p}_{a}\otimes \delta \tilde{p}_{b}
\end{equation*}%
\emph{is a (pseudo)metrical structure, then the real local functions: }%
\begin{equation*}
\begin{array}{l}
\left( \rho ,\eta \right) \overset{\ast }{H}_{\beta \gamma }^{\alpha
}=\left( \rho ,\eta \right) \overset{\ast }{\mathring{H}}_{\beta \gamma
}^{\alpha }+\displaystyle\frac{1}{2}\tilde{g}^{\alpha \varepsilon
}g_{\varepsilon \beta \overset{0}{\mid }\gamma },\vspace*{2mm} \\
\left( \rho ,\eta \right) \overset{\ast }{H}_{b\gamma }^{a}=\left( \rho
,\eta \right) \overset{\ast }{\mathring{H}}_{b\gamma }^{a}+\displaystyle%
\frac{1}{2}\tilde{g}_{be}g_{~\ \ \ \overset{0}{\mid }\gamma }^{ea},\vspace*{%
2mm} \\
\left( \rho ,\eta \right) \overset{\ast }{V}_{\beta }^{\alpha c}=\left( \rho
,\eta \right) \overset{\ast }{\mathring{V}}_{\beta }^{\alpha c}+\displaystyle%
\frac{1}{2}\tilde{g}^{\alpha \varepsilon }g_{\varepsilon \beta }\overset{0}{%
\mid }^{c},\vspace*{2mm} \\
\left( \rho ,\eta \right) \overset{\ast }{V}_{b}^{ac}=\left( \rho ,\eta
\right) \overset{\ast }{\mathring{V}}_{b}^{ac}+\displaystyle\frac{1}{2}%
\tilde{g}_{be}g^{ea}\overset{0}{\mid }^{c}%
\end{array}%
\leqno(5.8)
\end{equation*}%
\emph{are the components of a distinguished linear }$\left( \rho ,\eta
\right) $\emph{-connection such that the generalized tangent bundle }$\left(
\left( \rho ,\eta \right) T\overset{\ast }{E},\left( \rho ,\eta \right) \tau
_{\overset{\ast }{E}},\overset{\ast }{E}\right) $ \emph{becomes }$\left(
\rho ,\eta \right) $\emph{-(pseudo)metrizable.}

\textbf{Corollary 5.3 }\emph{In the particular case of Lie algebroids, }$%
\left( \eta ,h\right) =\left( Id_{M},Id_{M}\right) ,$\emph{\ then we obtain}%
\begin{equation*}
\begin{array}{l}
\rho \overset{\ast }{H}_{\beta \gamma }^{\alpha }=\rho \overset{\ast }{%
\mathring{H}}_{\beta \gamma }^{\alpha }+\displaystyle\frac{1}{2}\tilde{g}%
^{\alpha \varepsilon }g_{\varepsilon \beta \overset{0}{\mid }\gamma },%
\vspace*{2mm} \\
\rho \overset{\ast }{H}_{b\gamma }^{a}=\rho \overset{\ast }{\mathring{H}}%
_{b\gamma }^{a}+\displaystyle\frac{1}{2}\tilde{g}_{be}g_{~\ \ \ \overset{0}{%
\mid }\gamma }^{ea},\vspace*{2mm} \\
\rho \overset{\ast }{V}_{\beta }^{\alpha c}=\rho \overset{\ast }{\mathring{V}%
}_{\beta }^{\alpha c}+\displaystyle\frac{1}{2}\tilde{g}^{\alpha \varepsilon
}g_{\varepsilon \beta }\overset{0}{\mid }^{c},\vspace*{2mm} \\
\rho \overset{\ast }{V}_{b}^{ac}=\rho \overset{\ast }{\mathring{V}}_{b}^{ac}+%
\displaystyle\frac{1}{2}\tilde{g}_{be}g^{ea}\overset{0}{\mid }^{c}%
\end{array}%
\leqno(5.8)^{\prime }
\end{equation*}

\emph{In the classicale case, }$\left( \rho ,\eta ,h\right) =\left(
Id_{TE},Id_{M},Id_{M}\right) ,$\emph{\ then we obtain }%
\begin{equation*}
\begin{array}{l}
\overset{\ast }{H}_{jk}^{i}=\overset{\ast }{\mathring{H}}_{jk}^{i}+%
\displaystyle\frac{1}{2}\tilde{g}^{ih}g_{hj\overset{0}{\mid }k},\vspace*{2mm}
\\
\overset{\ast }{H}_{bk}^{a}=\overset{\ast }{\mathring{H}}_{bk}^{a}+%
\displaystyle\frac{1}{2}\tilde{g}_{be}g_{~\ \ \ \overset{0}{\mid }k}^{ea},%
\vspace*{2mm} \\
\overset{\ast }{V}_{j}^{ic}=\overset{\ast }{\mathring{V}}_{j}^{ic}+%
\displaystyle\frac{1}{2}\tilde{g}^{ih}g_{hj}\overset{0}{\mid }^{c},\vspace*{%
2mm} \\
\overset{\ast }{V}_{b}^{ac}=\overset{\ast }{\mathring{V}}_{b}^{ac}+%
\displaystyle\frac{1}{2}\tilde{g}_{be}g^{ea}\overset{0}{\mid }^{c}%
\end{array}%
\leqno(5.8)^{\prime \prime }
\end{equation*}

\section{Generalized Hamilton $\left( \protect\rho ,\protect\eta \right) $%
-spaces, Hamilton $\left( \protect\rho ,\protect\eta \right) $-spaces and
Cartan $\left( \protect\rho ,\protect\eta \right) $-spaces}

We consider the following diagram:
\begin{equation*}
\begin{array}{c}
\xymatrix{\overset{\ast }{E}\ar[d]_{\overset{\ast }{\pi }}&\left( F,\left[
,\right] _{F,h},\left( \rho ,\eta \right) \right)\ar[d]^\nu\\ M\ar[r]^h&N}%
\end{array}%
\end{equation*}%
such that $\left( E,\pi ,M\right) =\left( F,\nu ,N\right) $ and the
generalized tangent bundle
\begin{equation*}
\left( \left( \rho ,\eta \right) T\overset{\ast }{E},\left( \rho ,\eta
\right) \tau _{\overset{\ast }{E}},\overset{\ast }{E}\right)
\end{equation*}%
is $\left( \rho ,\eta \right) $-(pseudo)metrizable.

\bigskip\ Let%
\begin{equation*}
G=h_{ab}d\tilde{z}^{a}\otimes d\tilde{z}^{b}+g^{ab}\delta \tilde{p}%
_{a}\otimes \delta \tilde{p}_{b}
\end{equation*}%
be a (pseudo)metrical structure and
\begin{equation*}
\left( \left( \rho ,\eta \right) \overset{\ast }{H},\left( \rho ,\eta
\right) \overset{\ast }{V}\right)
\end{equation*}%
a distinguished linear $\left( \rho ,\eta \right) $-connection such that%
\emph{\ }%
\begin{equation*}
\begin{array}{c}
\left( \rho ,\eta \right) \overset{\ast }{D}_{X}G=0,~\forall X\in \Gamma
\left( \left( \rho ,\eta \right) T\overset{\ast }{E},\left( \rho ,\eta
\right) \tau _{\overset{\ast }{E}},\overset{\ast }{E}\right) .%
\end{array}%
\end{equation*}

\textbf{Definition 6.1 }A smooth \emph{Hamilton fundamental function} on the
dual vector bundle\break $\left( \overset{\ast }{E},\overset{\ast }{\pi }%
,M\right) $ is a mapping $\overset{\ast }{E}~\ ^{\underrightarrow{\ \ H\ \ }%
}~\ \mathbb{R}$ which satisfies the following conditions:\medskip

1. $H\circ \overset{\ast }{u}\in C^{\infty }\left( M\right) $, for any $%
\overset{\ast }{u}\in \Gamma \left( \overset{\ast }{E},\overset{\ast }{\pi }%
,M\right) \setminus \left\{ \overset{\ast }{0}\right\} $;\smallskip

2. $H\circ \overset{\ast }{0}\in C^{0}\left( M\right) $, where $\overset{%
\ast }{0}$ means the null section of $\left( \overset{\ast }{E},\overset{%
\ast }{\pi },M\right) .$\medskip

If $\left( U,\overset{\ast }{s}_{U}\right) $ is a local vector $\left(
m+r\right) $-chart for $\left( \overset{\ast }{E},\overset{\ast }{\pi }%
,M\right) $, then real function%
\begin{equation*}
\begin{array}[b]{c}
H^{ab}\overset{put}{=}\frac{\partial ^{2}H}{\partial p_{a}\partial p_{b}}%
\overset{put}{=}\frac{\partial }{\partial p_{a}}\left( \frac{\partial }{%
\partial p_{b}}\left( H\right) \right)%
\end{array}%
\end{equation*}%
is defined on $\overset{\ast }{\pi }^{-1}\left( U\right) $.

\textbf{Definition 6.2 }If for any local vector $m+r$-chart $\left( U,%
\overset{\ast }{s}_{U}\right) $ of $\left( \overset{\ast }{E},\overset{\ast }%
{\pi },M\right) ,$ we have:
\begin{equation*}
\begin{array}{c}
rank\left\Vert H^{ab}\left( \overset{\ast }{u}_{x}\right) \right\Vert =r,%
\end{array}%
\leqno(6.2)
\end{equation*}%
for any $\overset{\ast }{u}_{x}\in \overset{\ast }{\pi }^{-1}\left( U\right)
\backslash \left\{ \overset{\ast }{0}_{x}\right\} $, then we will say that
\emph{the Hamiltonian }$H$\emph{\ is regular.}

\textbf{Proposition 6.1} \emph{If the Hamiltonian }$H$\emph{\ is regular,
then for any local vector }$m+r$\emph{-chart }$\left( U,\overset{\ast }{s}%
_{U}\right) $ of $\left( \overset{\ast }{E},\overset{\ast }{\pi },M\right) ,$%
\emph{\ we obtain the real functions }$\tilde{H}_{ba}$\emph{\ locally
defined by}%
\begin{equation*}
\begin{array}{ccc}
\overset{\ast }{\pi }^{-1}\left( U\right) & ^{\underrightarrow{\ \ \tilde{H}%
_{ba}\ \ }} & \mathbb{R} \\
\overset{\ast }{u}_{x} & \longmapsto & \tilde{H}_{ba}\left( \overset{\ast }{u%
}_{x}\right)%
\end{array}%
,\leqno(6.3)
\end{equation*}%
\emph{where }$\left\Vert \tilde{H}_{ba}\left( \overset{\ast }{u}_{x}\right)
\right\Vert =\left\Vert H^{ab}\left( \overset{\ast }{u}_{x}\right)
\right\Vert ^{-1}$\emph{, for any }$\overset{\ast }{u}_{x}\in \overset{\ast }%
{\pi }^{-1}\left( U\right) \backslash \left\{ \overset{\ast }{0}_{x}\right\}
.$

\textbf{Definition 6.3 }A smooth \emph{Cartan fundamental function} on the
vector bundle $\left( \overset{\ast }{E},\overset{\ast }{\pi },M\right) $ is
a smooth Lagrange fundamental function $\overset{\ast }{E}~\ ^{%
\underrightarrow{\ \ K\ \ }}~\ \mathbb{R}_{+}$ which satisfies the following
conditions:\medskip

1. $K$ is positively $1$-homogenous on the fibres of vector bundle $\left(
\overset{\ast }{E},\overset{\ast }{\pi },M\right) ;$\smallskip

2. For any local vector $m+r$-chart $\left( U,\overset{\ast }{s}_{U}\right) $
of $\left( \overset{\ast }{E},\overset{\ast }{\pi },M\right) ,$ the hessian:%
\begin{equation*}
\left\Vert K^{2ab}\left( \overset{\ast }{u}_{x}\right) \right\Vert \leqno%
(6.4)
\end{equation*}%
is positively define for any $u_{x}\in \overset{\ast }{\pi }^{-1}\left(
U\right) \backslash \left\{ \overset{\ast }{0}_{x}\right\} $.

\textbf{Definition 6.4 }If the (pseudo)metrical structure $G$\ is determined
by a (pseudo)metrical structure
\begin{equation*}
\begin{array}{c}
g=g^{ab}d\tilde{p}_{a}\otimes d\tilde{p}_{b}\in \mathcal{T}~_{2}^{0}\left(
V\left( \rho ,\eta \right) T\overset{\ast }{E},\left( \rho ,\eta \right)
\tau _{\overset{\ast }{E}},\overset{\ast }{E}\right) ,%
\end{array}%
\end{equation*}%
namely
\begin{equation*}
G=\tilde{g}_{ab}d\tilde{z}^{a}\otimes d\tilde{z}^{b}+g^{ab}\delta \tilde{p}%
_{a}\otimes \delta \tilde{p}_{b},
\end{equation*}%
then the $\left( \rho ,\eta \right) $-(pseudo)metrizable vector bundle
\begin{equation*}
\begin{array}{c}
\left( \left( \rho ,\eta \right) T\overset{\ast }{E},\left( \rho ,\eta
\right) \tau _{\overset{\ast }{E}},\overset{\ast }{E}\right)%
\end{array}%
\end{equation*}%
will be called the \emph{generalized Hamilton }$\left( \rho ,\eta \right) $%
\emph{-space.}

In particular, if the (pseudo)metrical structure $g$\ is determined with the
help of a regular Hamilton (Cartan) fundamental function, namely $g=H^{ab}d%
\tilde{p}_{a}\otimes d\tilde{p}_{b}$ $\left( g=K^{2ab}d\tilde{p}_{a}\otimes d%
\tilde{p}_{b}\right) $, then the $\left( \rho ,\eta \right) $%
-(pseudo)metrizable vector bundle
\begin{equation*}
\left( \left( \rho ,\eta \right) T\overset{\ast }{E},\left( \rho ,\eta
\right) \tau _{\overset{\ast }{E}},\overset{\ast }{E}\right)
\end{equation*}%
will be called the \emph{Hamilton (Cartan) }$\left( \rho ,\eta \right) $%
\emph{-space.}

The generalized Hamilton $\left( Id_{T^{\ast }M},Id_{M}\right) $-spaces, the
Hamilton $\left( Id_{T^{\ast }M},Id_{M}\right) $-spaces, and the Cartan $%
\left( Id_{T^{\ast }M},Id_{M}\right) $-spaces are the usual generalized
Hamilton spaces, Hamilton spaces and Cartan spaces.

\textbf{Theorem 6.1 }\emph{If the (pseudo)metrical structure }$G$\emph{\ is
determined by a (pseudo)metrical structure }%
\begin{equation*}
\begin{array}{c}
g\in \mathcal{T}~_{2}^{0}\left( V\left( \rho ,\eta \right) T\overset{\ast }{E%
},\left( \rho ,\eta \right) \tau _{\overset{\ast }{E}},\overset{\ast }{E}%
\right) ,%
\end{array}%
\end{equation*}%
\emph{then, the real local functions:}
\begin{equation*}
(6.5)%
\begin{array}{ll}
\left( \rho ,\eta \right) \overset{\ast }{H}_{bc}^{a}\!\!\! & =\displaystyle%
\frac{1}{2}g^{ae}\left( \Gamma \left( \overset{\ast }{\tilde{\rho}},Id_{%
\overset{\ast }{E}}\right) \left( \overset{\ast }{\tilde{\delta}}_{b}\right)
\tilde{g}_{ec}+\Gamma \left( \overset{\ast }{\tilde{\rho}},Id_{\overset{\ast
}{E}}\right) \left( \overset{\ast }{\tilde{\delta}}_{c}\right) \tilde{g}%
_{be}\right. -\Gamma \left( \overset{\ast }{\tilde{\rho}},Id_{\overset{\ast }%
{E}}\right) \left( \overset{\ast }{\tilde{\delta}}_{e}\right) \tilde{g}_{bc}%
\vspace*{1mm} \\
& -\,\tilde{g}_{cd}L_{be}^{d}{\circ }h{\circ }\overset{\ast }{\pi }\left. +%
\tilde{g}_{bd}L_{ec}^{d}{\circ }h{\circ }\overset{\ast }{\pi }-\tilde{g}%
_{ed}L_{bc}^{d}{\circ }h{\circ }\overset{\ast }{\pi }\right) ,\vspace*{2mm}
\\
\left( \rho ,\eta \right) \overset{\ast }{V}_{b}^{ac}\!\!\! & =\displaystyle%
\frac{1}{2}\tilde{g}_{be}\left( \Gamma \left( \overset{\ast }{\tilde{\rho}}%
,Id_{\overset{\ast }{E}}\right) \left( \overset{\cdot }{\tilde{\partial}}%
^{c}\right) g^{ea}+\Gamma \left( \overset{\ast }{\tilde{\rho}},Id_{\overset{%
\ast }{E}}\right) \left( \overset{\cdot }{\tilde{\partial}}^{a}\right)
g^{ec}-\Gamma \left( \overset{\ast }{\tilde{\rho}},Id_{\overset{\ast }{E}%
}\right) \left( \overset{\cdot }{\tilde{\partial}}^{e}\right) g^{ac}\right)%
\end{array}%
\end{equation*}%
\medskip \emph{are the components of a normal distinguished linear }$\left(
\rho ,\eta \right) $\emph{-connection with }$\left( \rho ,\eta \right) $%
\emph{-}$\mathcal{H}\left( \mathcal{HH}\right) $\emph{\ and }$\left( \rho
,\eta \right) $\emph{-}$\mathcal{V}\left( \mathcal{VV}\right) $\emph{\
torsions free such that the generalized tangent bundle}\break $\left( \left(
\rho ,\eta \right) T\overset{\ast }{E},\left( \rho ,\eta \right) \tau _{%
\overset{\ast }{E}},\overset{\ast }{E}\right) $\emph{\ becomes generalized
Hamilton }$\left( \rho ,\eta \right) $\emph{-space.}\medskip

This normal distinguished linear $(\rho ,\eta )$-connection will be called
the \emph{generalized linear} $(\rho ,\eta )$\emph{-connection of
Levi-Civita type. }

\textbf{Corolary 6.1 }\emph{In the particular case of Lie algebroids, }$%
\left( \eta ,h\right) =\left( Id_{M},Id_{M}\right) ,$\emph{\ then we obtain}
\begin{equation*}
(6.5)^{\prime }%
\begin{array}{ll}
\rho \overset{\ast }{H}_{bc}^{a}\!\!\! & =\displaystyle\frac{1}{2}%
g^{ae}\left( \Gamma \left( \overset{\ast }{\tilde{\rho}},Id_{\overset{\ast }{%
E}}\right) \left( \overset{\ast }{\tilde{\delta}}_{b}\right) \tilde{g}%
_{ec}+\Gamma \left( \overset{\ast }{\tilde{\rho}},Id_{\overset{\ast }{E}%
}\right) \left( \overset{\ast }{\tilde{\delta}}_{c}\right) \tilde{g}%
_{be}\right. -\Gamma \left( \overset{\ast }{\tilde{\rho}},Id_{\overset{\ast }%
{E}}\right) \left( \overset{\ast }{\tilde{\delta}}_{e}\right) \tilde{g}_{bc}%
\vspace*{1mm} \\
& -\,\tilde{g}_{cd}L_{be}^{d}{\circ }\overset{\ast }{\pi }\left. +\tilde{g}%
_{bd}L_{ec}^{d}{\circ }\overset{\ast }{\pi }-\tilde{g}_{ed}L_{bc}^{d}{\circ }%
\overset{\ast }{\pi }\right) ,\vspace*{2mm} \\
\rho \overset{\ast }{V}_{b}^{ac}\!\!\! & =\displaystyle\frac{1}{2}\tilde{g}%
_{be}\left( \Gamma \left( \overset{\ast }{\tilde{\rho}},Id_{\overset{\ast }{E%
}}\right) \left( \overset{\cdot }{\tilde{\partial}}^{c}\right) g^{ea}+\Gamma
\left( \overset{\ast }{\tilde{\rho}},Id_{\overset{\ast }{E}}\right) \left(
\overset{\cdot }{\tilde{\partial}}^{a}\right) g^{ec}-\Gamma \left( \overset{%
\ast }{\tilde{\rho}},Id_{\overset{\ast }{E}}\right) \left( \overset{\cdot }{%
\tilde{\partial}}^{e}\right) g^{ac}\right)%
\end{array}%
\end{equation*}

\emph{In the classicale case, }$\left( \rho ,\eta ,h\right) =\left(
Id_{TE},Id_{M},Id_{M}\right) ,$\emph{\ then we obtain}
\begin{equation*}
\begin{array}{ll}
\overset{\ast }{H}_{bc}^{a}\!\!\! & =\frac{1}{2}g^{ae}\left( \overset{\ast }{%
\delta }_{b}\tilde{g}_{ec}+\overset{\ast }{\delta }_{c}\tilde{g}_{be}-%
\overset{\ast }{\delta }_{e}\tilde{g}_{bc}\vspace*{1mm}\right) \\
\overset{\ast }{V}_{b}^{ac}\!\!\! & =\displaystyle\frac{1}{2}\tilde{g}%
_{be}\left( \dot{\partial}^{c}g^{ea}+\dot{\partial}^{a}g^{ec}-\dot{\partial}%
^{e}g^{ac}\right)%
\end{array}%
\leqno(6.5)^{\prime \prime }
\end{equation*}

\emph{Moreover, if }$\left( E,\pi ,M\right) =\left( TM,\tau _{M},M\right) ,$%
\emph{\ then we obtain}%
\begin{equation*}
\begin{array}{ll}
\overset{\ast }{H}_{jk}^{i}\!\!\! & =\frac{1}{2}g^{ih}\left( \overset{\ast }{%
\delta }_{j}\tilde{g}_{hk}+\overset{\ast }{\delta }_{k}\tilde{g}_{jh}-%
\overset{\ast }{\delta }_{h}\tilde{g}_{jk}\vspace*{1mm}\right) \\
\overset{\ast }{V}_{j}^{ik}\!\!\! & =\displaystyle\frac{1}{2}\tilde{g}%
_{jh}\left( \dot{\partial}^{k}g^{hi}+\dot{\partial}^{i}g^{hk}-\dot{\partial}%
^{h}g^{ik}\right)%
\end{array}%
\leqno(6.5)^{\prime \prime \prime }
\end{equation*}

\textbf{Theorem 6.2 }\emph{Let\ }$\left( \left( \rho ,\eta \right) \overset{%
\ast }{H},\left( \rho ,\eta \right) \overset{\ast }{V}\right) $\emph{\ be
the normal dis\-tin\-guished linear }$\left( \rho ,\eta \right) $\emph{%
-connec\-tion presented in the previous theorem. If }%
\begin{equation*}
\overset{\ast }{\mathbb{T}}_{bc}^{a}\tilde{\delta}_{a}\otimes d\tilde{z}%
^{b}\otimes d\tilde{z}^{c}\in \mathcal{T}_{20}^{10}\left( \left( \rho ,\eta
\right) T\overset{\ast }{E},\left( \rho ,\eta \right) \tau _{\overset{\ast }{%
E}},\overset{\ast }{E}\right)
\end{equation*}%
\emph{and }%
\begin{equation*}
\overset{\ast }{\mathbb{S}}_{b}^{ac}\overset{\cdot }{\tilde{\partial}}%
^{b}\otimes \delta \tilde{p}_{a}\otimes \delta \tilde{p}_{c}\in \mathcal{T}%
_{01}^{02}\left( \left( \rho ,\eta \right) T\overset{\ast }{E},\left( \rho
,\eta \right) \tau _{\overset{\ast }{E}},\overset{\ast }{E}\right)
\end{equation*}%
\emph{such that they satisfy the conditions:}%
\begin{equation*}
\overset{\ast }{\mathbb{T}}_{bc}^{a}=-\overset{\ast }{\mathbb{T}}_{cb}^{a},~%
\overset{\ast }{\mathbb{S}}_{b}^{ac}=-\overset{\ast }{\mathbb{S}}%
_{b}^{ac},~\forall ab,c\in \overline{1,r},
\end{equation*}%
\emph{then the following real local functions:\ }%
\begin{equation*}
\begin{array}{l}
\left( \rho ,\eta \right) \overset{\ast }{\tilde{H}}_{bc}^{a}=\left( \rho
,\eta \right) \overset{\ast }{H}_{bc}^{a}+\displaystyle\frac{1}{2}%
g^{ae}\left( \tilde{g}_{ed}\overset{\ast }{\mathbb{T}}_{bc}^{d}-\tilde{g}%
_{bd}\overset{\ast }{\mathbb{T}}_{ec}^{d}+\tilde{g}_{cd}\overset{\ast }{%
\mathbb{T}}_{be}^{d}\right) ,\vspace*{2mm} \\
\left( \rho ,\eta \right) \overset{\ast }{\tilde{V}}_{b}^{ac}=\left( \rho
,\eta \right) \overset{\ast }{V}_{b}^{ac}+\displaystyle\frac{1}{2}\tilde{g}%
_{be}\cdot \left( g^{ed}\overset{\ast }{\mathbb{S}}_{d}^{ac}-g^{ad}\overset{%
\ast }{\mathbb{S}}_{d}^{ec}+g^{cd}\overset{\ast }{\mathbb{S}}%
_{d}^{ae}\right)
\end{array}%
\leqno(6.6)
\end{equation*}%
\emph{are the components of a normal distinguished linear }$\left( \rho
,\eta \right) $\emph{-connection with }$\left( \rho ,\eta \right) $\emph{-}$%
\mathcal{H}\left( \mathcal{HH}\right) $\emph{\ and }$\left( \rho ,\eta
\right) $\emph{-}$\mathcal{V}\left( \mathcal{VV}\right) $\emph{\ torsions a
priori given such that the generalized tangent bundle}\break $\left( \left(
\rho ,\eta \right) T\overset{\ast }{E},\left( \rho ,\eta \right) \tau _{%
\overset{\ast }{E}},\overset{\ast }{E}\right) $\emph{\ becomes generalized
Hamilton }$\left( \rho ,\eta \right) $\emph{-space.}

\emph{Moreover, we obtain: }%
\begin{equation*}
\begin{array}{l}
\overset{\ast }{\mathbb{T}}_{bc}^{a}=\left( \rho ,\eta \right) \overset{\ast
}{\tilde{H}}_{bc}^{a}-\left( \rho ,\eta \right) \overset{\ast }{\tilde{H}}%
_{cb}^{a}-L_{bc}^{a}\circ h\circ \overset{\ast }{\pi },\vspace*{2mm} \\
\overset{\ast }{\mathbb{S}}_{b}^{ac}=\left( \rho ,\eta \right) \overset{\ast
}{\tilde{V}}_{b}^{ac}-\left( \rho ,\eta \right) \overset{\ast }{\tilde{V}}%
_{b}^{ca}.%
\end{array}%
\leqno\left( 6.7\right)
\end{equation*}

\textbf{Corollary 6.2} \emph{In the particular case of Lie algebroids, }$%
\left( \eta ,h\right) =\left( Id_{M},Id_{M}\right) ,$\emph{\ then we obtain}%
\begin{equation*}
\begin{array}{l}
\rho \overset{\ast }{\tilde{H}}_{bc}^{a}=\rho \overset{\ast }{H}_{bc}^{a}+%
\displaystyle\frac{1}{2}g^{ae}\left( \tilde{g}_{ed}\overset{\ast }{\mathbb{T}%
}_{bc}^{d}-\tilde{g}_{bd}\overset{\ast }{\mathbb{T}}_{ec}^{d}+\tilde{g}_{cd}%
\overset{\ast }{\mathbb{T}}_{be}^{d}\right) ,\vspace*{2mm} \\
\rho \overset{\ast }{\tilde{V}}_{b}^{ac}=\rho \overset{\ast }{V}_{b}^{ac}+%
\displaystyle\frac{1}{2}\tilde{g}_{be}\cdot \left( g^{ed}\overset{\ast }{%
\mathbb{S}}_{d}^{ac}-g^{ad}\overset{\ast }{\mathbb{S}}_{d}^{ec}+g^{cd}%
\overset{\ast }{\mathbb{S}}_{d}^{ae}\right) .%
\end{array}%
\leqno(6.6)^{\prime }
\end{equation*}%
\emph{and }%
\begin{equation*}
\begin{array}{l}
\overset{\ast }{\mathbb{T}}_{bc}^{a}=\rho \overset{\ast }{\tilde{H}}%
_{bc}^{a}-\rho \overset{\ast }{\tilde{H}}_{cb}^{a}-L_{bc}^{a}\circ \overset{%
\ast }{\pi },\vspace*{2mm} \\
\overset{\ast }{\mathbb{S}}_{b}^{ac}=\rho \overset{\ast }{\tilde{V}}%
_{b}^{ac}-\rho \overset{\ast }{\tilde{V}}_{b}^{ca}.%
\end{array}%
\leqno\left( 6.7\right) ^{\prime }
\end{equation*}

\emph{In the classicale case, }$\left( \rho ,\eta ,h\right) =\left(
Id_{TE},Id_{M},Id_{M}\right) ,$\emph{\ then we obtain}
\begin{equation*}
\begin{array}{l}
\overset{\ast }{\tilde{H}}_{bc}^{a}=\overset{\ast }{H}_{bc}^{a}+\displaystyle%
\frac{1}{2}g^{ae}\left( \tilde{g}_{ed}\overset{\ast }{\mathbb{T}}_{bc}^{d}-%
\tilde{g}_{bd}\overset{\ast }{\mathbb{T}}_{ec}^{d}+\tilde{g}_{cd}\overset{%
\ast }{\mathbb{T}}_{be}^{d}\right) ,\vspace*{2mm} \\
\overset{\ast }{\tilde{V}}_{b}^{ac}=\overset{\ast }{V}_{b}^{ac}+\displaystyle%
\frac{1}{2}\tilde{g}_{be}\cdot \left( g^{ed}\overset{\ast }{\mathbb{S}}%
_{d}^{ac}-g^{ad}\overset{\ast }{\mathbb{S}}_{d}^{ec}+g^{cd}\overset{\ast }{%
\mathbb{S}}_{d}^{ae}\right) .%
\end{array}%
\leqno(6.6)^{\prime \prime }
\end{equation*}%
\emph{and }%
\begin{equation*}
\begin{array}{l}
\overset{\ast }{\mathbb{T}}_{bc}^{a}=\overset{\ast }{\tilde{H}}_{bc}^{a}-%
\overset{\ast }{\tilde{H}}_{cb}^{a},\vspace*{2mm} \\
\overset{\ast }{\mathbb{S}}_{b}^{ac}=\overset{\ast }{\tilde{V}}_{b}^{ac}-%
\overset{\ast }{\tilde{V}}_{b}^{ca}.%
\end{array}%
\leqno\left( 6.7\right) ^{\prime \prime }
\end{equation*}

\emph{In particular, if }$\left( E,\pi ,M\right) =\left( TM,\tau
_{M},M\right) ,$\emph{\ then we obtain}
\begin{equation*}
\begin{array}{l}
\overset{\ast }{\tilde{H}}_{jk}^{i}=\overset{\ast }{H}_{jk}^{i}+\displaystyle%
\frac{1}{2}g^{ih}\left( \tilde{g}_{hl}\overset{\ast }{\mathbb{T}}_{jk}^{l}-%
\tilde{g}_{jl}\overset{\ast }{\mathbb{T}}_{hk}^{l}+\tilde{g}_{kl}\overset{%
\ast }{\mathbb{T}}_{jh}^{l}\right) ,\vspace*{2mm} \\
\overset{\ast }{\tilde{V}}_{j}^{ik}=\overset{\ast }{V}_{j}^{ik}+\displaystyle%
\frac{1}{2}\tilde{g}_{jh}\cdot \left( g^{hl}\overset{\ast }{\mathbb{S}}%
_{l}^{ik}-g^{il}\overset{\ast }{\mathbb{S}}_{l}^{hk}+g^{kl}\overset{\ast }{%
\mathbb{S}}_{l}^{ih}\right) .%
\end{array}%
\leqno(6.6)^{\prime \prime \prime }
\end{equation*}%
\emph{and }%
\begin{equation*}
\begin{array}{l}
\overset{\ast }{\mathbb{T}}_{jk}^{i}=\overset{\ast }{\tilde{H}}_{jk}^{i}-%
\overset{\ast }{\tilde{H}}_{kj}^{i},\vspace*{2mm} \\
\overset{\ast }{\mathbb{S}}_{j}^{ik}=\overset{\ast }{\tilde{V}}_{j}^{ik}-%
\overset{\ast }{\tilde{V}}_{j}^{ki}.%
\end{array}%
\leqno\left( 6.7\right) ^{\prime \prime \prime }
\end{equation*}

\hfill
\begin{tabular}{c}
SECONDARY SCHOOL \textquotedblleft CORNELIUS RADU\textquotedblright , \\
RADINESTI VILLAGE, 217196, GORJ COUNTY, ROMANIA \\
e-mail: c\_arcus@yahoo.com, c\_arcus@radinesti.ro%
\end{tabular}

\end{document}